\documentclass[leqno,11pt]{article}

\usepackage{amssymb,amsmath,amscd,amsfonts,a4,psfrag,graphicx,latexsym,epsfig,color}
\usepackage[all]{xy}

\usepackage{tikz}

\usepackage{verbatim}

\usetikzlibrary{matrix}

\newtheorem{theorem}{Theorem}[section]
\newtheorem{lemma}[theorem]{Lemma}
\newtheorem{proposition}[theorem]{Proposition}
\newtheorem{corollary}[theorem]{Corollary}
\newtheorem{definition}[theorem]{Definition}
\newtheorem{remark}[theorem]{Remark}

\newcommand{\agot}{\ensuremath{\mathfrak{a}}}

\newcommand{\ggot}{\ensuremath{\mathfrak{g}}}

\newcommand{\kgot}{\ensuremath{\mathfrak{k}}}

\newcommand{\pgot}{\ensuremath{\mathfrak{p}}}
\newcommand{\qgot}{\ensuremath{\mathfrak{q}}}
\newcommand{\ngot}{\ensuremath{\mathfrak{n}}}

\newcommand{\sgot}{\ensuremath{\mathfrak{s}}}
\newcommand{\tgot}{\ensuremath{\mathfrak{t}}}
\newcommand{\ugot}{\ensuremath{\mathfrak{u}}}
\newcommand{\Rgot}{\ensuremath{\mathfrak{R}}}

\newcommand{\Bcal}{\ensuremath{\mathcal{B}}}
\newcommand{\Ccal}{\ensuremath{\mathcal{C}}}
\newcommand{\Dcal}{\ensuremath{\mathcal{D}}}

\newcommand{\Fcal}{\ensuremath{\mathcal{F}}}

\newcommand{\Lcal}{\ensuremath{\mathcal{L}}}

\newcommand{\Tcal}{\ensuremath{\mathcal{T}}}
\newcommand{\Xcal}{\ensuremath{\mathcal{X}}}
\newcommand{\Ycal}{\ensuremath{\mathcal{Y}}}
\newcommand{\Zcal}{\ensuremath{\mathcal{Z}}}
\newcommand{\Ucal}{\ensuremath{\mathcal{U}}}
\newcommand{\Vcal}{\ensuremath{\mathcal{V}}}

\newcommand{\J}{\ensuremath{\mathbb{J}}}
\newcommand{\Z}{\ensuremath{\mathbb{Z}}}
\newcommand{\C}{\ensuremath{\mathbb{C}}}

\newcommand{\Q}{\ensuremath{\mathbb{Q}}}
\newcommand{\R}{\ensuremath{\mathbb{R}}}

\newcommand{\T}{\ensuremath{\hbox{\bf T}}}

\newcommand{\grad}{{\ensuremath{\rm grad}}}
\newcommand{\critical}{{\ensuremath{\rm crit}}}

\newcommand{\mini}{\ensuremath{\hbox{\scriptsize \rm s}}}

\newcommand{\reg}{\operatorname{\hbox{\rm \tiny reg}}}

\newcommand{\infrp}{\operatorname{\hbox{\rm \tiny inf-RP}}}

\newcommand{\rp}{\operatorname{\hbox{\rm \tiny RP}}}

\begin{document}

\title{Moment polytopes in real symplectic geometry I}

\author{Paul-Emile Paradan\footnote{IMAG, Univ Montpellier, CNRS, email : paul-emile.paradan@umontpellier.fr}}

\maketitle

\date{}


\begin{abstract}
Let $Z$ be the real part of a K\"{a}hler Hamiltonian manifold $M$. The O'Shea-Sjamaar's Theorem tells us that the moment polytope
$\Delta_\pgot(Z)$ corresponds to the anti-invariant part of the Kirwan polytope $\Delta_\ugot(M)$. 
The purpose of the present paper is to explain how to parameterize the equations of the facets 
of $\Delta_\pgot(Z)$ in terms of {\em real Ressayre's pairs of $Z$}.

\end{abstract}



\tableofcontents

\section{Introduction}

Let $U_\C$ be a connected complex reductive group equipped with a complex-conjugate involution $\sigma$. Let 
$U$ be the maximal compact subgroup of $U_\C$ that is associated with a Cartan involution commuting to $\sigma$ : so 
$U$ is invariant under the involution and we denote by $K$ the connected component of the fixed point subgroup $U^\sigma$. At the level of Lie algebras we have a decomposition 
$\ugot=\ugot^{\sigma}\oplus\ugot^{-\sigma}$ where $\ugot^{\pm\sigma}=\{X\in\ugot, \sigma(X)=\pm X\}$. 

Let $G$ be the connected component of the closed subgroup of $U_\C$ fixed by $\sigma$ : it is a real reductive subgroup that is stable under the Cartan involution. 
Hence we have the isomorphism $K\times\pgot \to G$, $(k,X)\mapsto k e^X$ where $\pgot=i \ugot^{-\sigma}$.

We consider now a proper K\"{a}hler Hamiltonian $U$-manifold $(M,\Omega,\J)$. By that we mean~:
\begin{itemize}
\item  $U_\C$ acts holomorphically on the complex manifold $(M,\J)$,
\item the K\"{a}hler form $\Omega$ is $U$-invariant, 
\item there is a $U$-equivariant {\em proper} moment mapping $\Phi_\ugot : M\to \ugot^*$ satisfying the relation 
\begin{equation}\label{eq:hamiltonian}
d\langle\Phi_\ugot,X\rangle=-\Omega(X_M, -),\quad \forall X\in\ugot. 
\end{equation}
Here $X_M(m)=-X\cdot m =\frac{d}{dt}\vert_{t=0}e^{-tX}m$ is the vector field generated by $X\in \ugot$.
\end{itemize}

In this paper we will give a particular attention to the $K$-equivariant map 
$$
\Phi_\pgot : M\to \pgot^*
$$ 
defined by the relations
$\langle\Phi_\pgot,\beta\rangle=\langle\Phi_\ugot,i\beta\rangle$, $\forall \beta\in\pgot$. If we denote by $j:\pgot\to\ugot^{-\sigma}$ the isomorphism $X\mapsto iX$, we have $\Phi_\pgot=j^*\circ \pi\circ\Phi_{\ugot}$, where $j^*:(\ugot^{-\sigma})^*\to \pgot^*$ is the dual map and $\pi:\ugot^*\to(\ugot^{-\sigma})^*$ is the canonical projection. Thanks to (\ref{eq:hamiltonian}), we see that $\beta_M$ is the gradient vector 
field of the function $\langle\Phi_\pgot,\beta\rangle$ for any $\beta\in\pgot$.

We suppose now that the K\"{a}hler Hamiltonian $U$-manifold $(M,\Omega,\J)$ possesses a {\em real structure}, that is an anti-holomorphic involution $\tau :M\to M$ satisfying the following conditions : $\tau^*(\Omega)=-\Omega$ and 
\begin{itemize}
\item $\tau(g\cdot m)=\sigma(g)\cdot\tau(m)$,
\item $\Phi_\ugot(\tau(m))=-\sigma(\Phi_\ugot(m))$,
\end{itemize}
for any $(g,m)\in U^\C\times M$. 

Consider the submanifold $Z=\{m\in M, \tau(m)=m\}$ : if $Z$ is {\em non-empty}, it is a Lagrangian submanifold of $(M,\Omega)$, not necessarily connected, stable under the action of the real reductive group $G$. 
In this context, we have the following important result \cite{OSS99}.

\begin{theorem}[O'Shea-Sjamaar]\label{theo:OSS}
$$
\Phi_\ugot(M)\bigcap (\ugot^{-\sigma})^*=\Phi_\ugot(Z)\underset{j^*}{\simeq} \Phi_\pgot(Z).
$$
\end{theorem}

Let us rephrase this theorem in terms of moment polytopes. Let $\agot\subset \pgot$ be a \break maximal abelian subspace, and choose a $\sigma$-invariant maximal torus 
$T$ of $U$ with Lie algebra $\tgot$ containing $i\agot$. We can define Weyl chambers\footnote{$\ugot^*/U\simeq \tgot^*_+$ and $\pgot^*/K\simeq \agot^*_+$} $\tgot^*_+\subset \tgot^*$ 
and $\agot^*_+\subset \agot^*$ such that \break $\tgot^*_+\cap (\ugot^{-\sigma})^*{\simeq}_{j^*} \agot^*_+ $ (see the Appendix in \cite{OSS99}).

We consider the moment polytopes
$$
\Delta_\ugot(M)=\Phi_\ugot(M)\cap \tgot^*_+\qquad {\rm and}\qquad \Delta_\pgot(Z)=\Phi_\pgot(Z)\cap \agot^*_+.
$$
The O'Shea-Sjamaar's theorem  can be rephrased by saying that the intersection $\Delta_\ugot(M)\cap(\ugot^{-\sigma})^*$ coincides with $\Delta_\pgot(Z)$ through the linear isomorphism 
$j^*$.

The Convexity Theorem \cite{Atiyah82,Guillemin-Sternberg82.bis,Kirwan.84.bis,sjamaar98} tells us that $\Delta_\ugot (M) $ is a closed convex locally polyhedral subset. Thanks to Theorem \ref{theo:OSS}, we see that
$\Delta_\pgot (Z) $ is also  a closed convex locally polyhedral subset.  In a previous article \cite{pep-RP}, we explained how to parameterize the facets of 
$\Delta_\ugot(M) $ in terms of {\em Ressayre's pairs}, a notion that N. Ressayre introduces in the algebraic setting \cite{Ressayre-10}.

The goal of this article is to show that we can adapt this notion in the real setting in order to parameterize the facets of $\Delta_\pgot (Z)$. A future article \cite{pep-real-RP-2} 
will be devoted to the study of examples.

\medskip

Since the founding article of R. Richardson and P. J. Slodowy \cite{Richardson-Slodowy90} where they showed that the Kempf-Ness theorem extends to representations of real reductive 
groups, many authors have studied extensions of geometric invariant theory to the real framework \cite{Heinzner-Schwarz07,HS10,Miebach-Stotzel10,Bohm-Lafuente17,Biliotti20}.

Let us cite the work of J. Lauret \cite{Lauret10} and Heinzner-Schwarz-St\"{o}tzel \cite{HSS08} where they obtain a stratification theorem \`a la Kirwan-Ness for real 
reductive group actions on real manifolds. In the present work, the existence of real stratifications  will be the main tool to parameterize the facets of the moment 
polytope $\Delta_\pgot(Z)$.  As we are working in a less general framework than that of \cite{HSS08}, we propose another proof of their result 
based on the Kirwan-Ness stratification of the ambiant complex manifold (see \S \ref{sec:real-stratification}). 
This method will allow us to show that these real stratifications admit an open and dense stratum. 
This key point provides a simple proof of O'Shea-Sjamaar's theorem in the K\"{a}hler setting (see \S \ref{sec:refinement}).

\section{Statement of the main result}

For the rest of this section, we suppose that $Z\neq \emptyset$ and we fix a connected component $\Zcal$ of $Z$ : it is a Riemannian manifold equipped with an action of the groups 
$K$ and  $G$.

\subsection{A refinement of O'Shea-Sjamaar's Theorem} 

The following result is proved in \S \ref{sec:refinement}.
\begin{theorem}
The following relations hold :
\begin{itemize}
\item $\Delta_\ugot(M)\cap(\ugot^{-\sigma})^*\underset{j^*}{\simeq}\Delta_\pgot(Z)$,
\item $\Delta_\pgot(\Zcal)=\Delta_\pgot(Z)$.
\end{itemize}
\end{theorem}

The equality $\Delta_\pgot(\Zcal)=\Delta_\pgot(Z)$ was already obtained by L. O'Shea and R. Sjamaar when $M$ is a polarized algebraic variety (see Corollary 5.11 in \cite{OSS99}).

\subsection{Admissible elements} 

We start by introducing the notion of {\em admissible elements}. The group \break ${\rm Hom}(U(1),T)$ admits a natural identification with the lattice $\wedge:=\frac{1}{2\pi}\ker(\exp : \tgot\to T)$ of the vector space $\tgot$. A vector $\gamma\in \agot$ is called {\em rational} if $i\gamma$ belongs to the $\Q$-vector space $\tgot_\Q$ generated by $\wedge$.

The stabilizer subgroups of $m\in M$ relatively to the $K$ and $G$ actions are denoted respectively by $K_m$ and $G_m$ : their Lie algebras are denote by $\kgot_m$ and $\ggot_m$. We will take a particular attention to the subspace $\pgot_m=\{X\in\pgot, X\cdot m=0\}\subset\ggot_m$. 
 If $\gamma\in\agot$, we denote by $\Zcal^\gamma$ the submanifold where the vector field $z\mapsto \gamma\cdot z$ vanishes.

\begin{definition} Let us define 
$$
\dim_\pgot(\Xcal):=\min_{z\in\Xcal}\,\dim(\pgot_z)
$$
for any subset $\Xcal\subset \Zcal$. A non-zero element $\gamma\in\agot$ is called {\em admissible} if $\gamma$ is rational, and if 
$\dim_\pgot(\Zcal^\gamma)-\dim_\pgot(\Zcal)\in\{0,1\}$.
\end{definition}

\subsection{Real Ressayre's pair}

The aim of this section is to introduce the notion of real Ressayre's pair on the $G$-manifold $\Zcal$.

For any $z\in \Zcal$, the infinitesimal action of $\ggot$ on $\Zcal$ defines a real linear map 
\begin{eqnarray}\label{eq:rho-z}
\rho_z:\ggot & \longrightarrow & \T_z \Zcal\\
X & \longmapsto & X\cdot m \nonumber
\end{eqnarray}
that is equivariant under the action of the stabilizer subgroup $G_z$.

\begin{definition}
Consider a symmetric endomorphism $\Lcal(\gamma)$ of an euclidien vector space $E$. We associate the eigenspace
$E^{\gamma=a}=\{v\in E, \Lcal(\gamma)v= av\}$ to any $a\in\R$.  We have the decomposition
$E=E^{\gamma>0}\oplus E^{\gamma=0}\oplus E^{\gamma<0}$ where $E^{\gamma>0}=\sum_{a>0}E^{\gamma=a}$, and $E^{\gamma<0}=\sum_{a<0}E^{\gamma=a}$.
\end{definition}

\medskip

Let $\Sigma$ be the set of (non-zero) roots relative to the action of $\agot$ on $\ggot$. The choice of the Weyl chamber $\agot^*_+$ induces a set  
$\Sigma^+$ of positive roots, and thus a decomposition $\ggot=\ngot_{-}\oplus \kgot'\oplus\agot \oplus\ngot$ where 
$\ngot=\sum_{\alpha>0}\ggot_\alpha$, $\ngot_{-}=\sum_{\alpha<0}\ggot_\alpha$ and $\kgot'$ is the centralizer subalgebra of $\agot$ in $\kgot$.

\begin{definition}\label{def:parabolic-Q}
We denote by $Q\subset G$ the (minimal) parabolic subgroup with Lie algebra $\qgot:= \kgot'\oplus\agot \oplus\ngot$.
\end{definition}

Consider $(x,\gamma)\in \Zcal\times\agot$ such that $x\in \Zcal^\gamma$. The infinitesimal action of $\gamma$ defines 
symmetric endomorphisms $\Lcal(\gamma):\T_x \Zcal\to \T_x \Zcal$ and $\Lcal(\gamma):\ggot\to \ggot$ satisfying $\Lcal(\gamma)\circ\rho_x=\rho_x\circ\Lcal(\gamma)$.
Thus the morphism (\ref{eq:rho-z}) induces a real linear map 
\begin{equation}\label{eq:rho-gamma}
\rho_x^\gamma:\ngot^{\gamma>0}\longrightarrow  (\T_x \Zcal)^{\gamma>0}.
\end{equation}

\medskip

\begin{definition}\label{def:infinitesimal-real-ressayre-pair}
Let $\gamma\in\agot$ be a non-zero element, and let $\Ccal\subset \Zcal^\gamma$ be a connected component. The data $(\gamma,\Ccal)$ is called 
an {\bf infinitesimal real Ressayre's pair} of $\Zcal$ if $\exists x\in \Ccal$, such as $\rho_x^\gamma$ is an isomorphism. 
If furthermore we have  $\dim_\pgot(\Ccal)-\dim_\pgot(\Zcal)\in\{0,1\}$, and $\gamma$ is rational, we call $(\gamma,\Ccal)$ a {\em regular infinitesimal real Ressayre's pair} of $\Zcal$.
\end{definition}

\begin{remark}
When $U=T$ is abelian, a couple $(\gamma,\Ccal)$ is an infinitesimal real Ressayre's pair when the vector bundle $(\T\Zcal\vert_\Ccal)^{\gamma>0}$ is equal to the ``zero'' bundle.
\end{remark}

\medskip

Let us now introduce a more restrictive notion, that of real Ressayre's pair. 

Let $\gamma\in\agot$ be a non-zero element. The fonction $\langle\Phi_\pgot,\gamma\rangle: \Zcal^\gamma\to\R$ is locally constant. Let 
$\Ccal=\Ccal_1\cup\cdots\cup\Ccal_p\subset \Zcal^\gamma$ be a {\em union} of connected components such that $\langle\Phi_\pgot,\gamma\rangle$ is constant on $\Ccal$: we denote 
by $\langle\Phi_\pgot(\Ccal),\gamma\rangle$ its value. We consider the {\em real} Bialynicki-Birula's submanifold\footnote{See \S \ref{sec:BB}.} 
\begin{equation}\label{eq:BB}
\Ccal^-:=\{z\in \Zcal,  \lim_{t\to\infty} \exp(t\gamma) z\ \in \Ccal\}.
\end{equation}
We see that for any $x\in \Ccal$, $(\T_x \Zcal)^{\gamma\leq 0}=\T_x \Ccal^-$. 
Consider now the real parabolic subgroup $P_\gamma\subset G$ defined by
\begin{equation}\label{eq:P-gamma-reel}
P_\gamma=\{g\in G, \lim_{t\to\infty}\exp(t\gamma)g\exp(-t\gamma)\ {\rm exists}\}.
\end{equation}
The submanifold $\Ccal^-$ is invariant under the action of $P_\gamma$, hence we can  consider the manifold 
$Q\times_{Q\cap P_\gamma} \Ccal^-$ and the smooth map 
$$
{\rm q}_\gamma: Q\times_{Q\cap P_\gamma} \Ccal^-\to \Zcal
$$
that sends $[q,x]$ to $qx$. We immediately see that for any $x\in \Ccal$, the tangent map 
$\T{\rm q}_\gamma\vert_x$ is an isomorphism if and only if $\rho_x^\gamma$ is an isomorphism.

\begin{definition}\label{def:real-ressayre-pair}
Let $\gamma\in\agot$ be a non-zero element, and let $\Ccal\subset \Zcal^\gamma$ be a  union of connected components such that 
$\langle\Phi_\pgot,\gamma\rangle$ is constant on $\Ccal$. The data $(\gamma,\Ccal)$ is called a {\bf real Ressayre's pair} of $\Zcal$ 
if the following conditions hold
\begin{itemize}
\item The image of ${\rm q}_{\gamma}$ contains a dense open subset of $\Zcal$.
\item There exists a $Q\cap P_{\gamma}$-invariant, open and dense subset $U\subset \Ccal^-$, intersecting $\Ccal$, 
so that ${\rm q}_{\gamma}$ defines a diffeomorphism $Q\times_{Q\cap P_{\gamma}}U\simeq Q U$.
\end{itemize}
If furthermore we have  $\dim_\pgot(\Ccal)-\dim_\pgot(\Zcal)\in\{0,1\}$, and $\gamma$ is rational, we call $(\gamma,\Ccal)$ a regular real Ressayre's pair.
\end{definition}

\begin{remark}
When $U=T$ is abelian, a couple $(\gamma,\Ccal)$ is a real Ressayre's pair when the real Bialynicki-Birula's submanifold $\Ccal^-$ is open and dense in $\Zcal$.
\end{remark}

\subsection{Main result}

The main result of this article is the following theorem.

\begin{theorem}\label{th:real-ressayre-pairs} 
For $\xi\in\agot^*_{+}$, the following statements are equivalent:
\begin{enumerate}
\item $\xi\in\Delta_\pgot(Z)$.
\item $\xi\in\Delta_\pgot(\Zcal)$.
\item For any regular infinitesimal real Ressayre's pair $(\gamma,\Ccal)$ of $\Zcal$, we have $\langle \xi,\gamma\rangle\geq \langle \Phi_\pgot(\Ccal),\gamma\rangle$.
\item For any regular real Ressayre's pair $(\gamma,\Ccal)$ of $\Zcal$, we have $\langle \xi,\gamma\rangle\geq \langle \Phi_\pgot(\Ccal),\gamma\rangle$.
\end{enumerate}
\end{theorem}

\begin{remark}
In the previous theorem, the result still holds if we drop the ``regular" hypothesis on $(\gamma,\Ccal)$.
\end{remark}

\section{Kirwan-Ness stratifications : complex and real settings}

Let us now choose a rational $U$-invariant inner product on $\ugot_\C$ that is invariant under the involution $\sigma$ (see the Appendix).
By rational we mean that for the maximal torus $T\subset U$ with Lie algebra $\tgot$, 
the inner product takes integral values on the lattice $\wedge:=\frac{1}{2\pi}\ker(\exp:\tgot\to T)$. Let us denote by $\wedge^*\subset \tgot^*$ the dual lattice : 
$\wedge^*=\hom(\wedge,\Z)$. We associate to the lattices $\wedge$ and $\wedge^*$ the $\Q$-vector space $\tgot_\Q$ and $\tgot^*_\Q$ generated by them: 
the vectors belonging to $\tgot_\Q$ and $\tgot^*_\Q$ are designed as rational.

The invariant scalar product on $\ugot$ induces an identification $\ugot^*\simeq \ugot,\xi\mapsto \xi^\flat $ such as $\tgot_\Q\simeq\tgot^*_\Q$. 
To simplify our notation, we will not distinguish between $\xi$ and $\xi^\flat$: for example we write $M^{\lambda}$ for the submanifold fixed by 
the subgroup generated by $\lambda^\flat$.

We come back to the setting of a  K\"{a}hler Hamiltonian $U$-manifold $(M,\Omega,\J)$ with proper moment map $\Phi_\ugot:M\to\ugot^*$ and an anti-holomorphic involution $\tau$.
We fix a connected component $\Zcal$ of  the real part $Z=M^\tau$.

\subsection{Bialynicki-Birula's submanifolds}\label{sec:BB}

Let us consider an element $\lambda\in\ugot$ and a connected component
$C$ of the complex submanifold $M^\lambda:=\{m\in M, \gamma\cdot m=0\}$. As in the introduction, we define the 
subset $C^-:=\{m\in M,  \lim_{t\to\infty} \exp(-it\lambda)m\ \in C\}$ and the projection $p_C:C^-\to C$ that sends $m\in C^-$ to 
$\lim_{t\to\infty} \exp(-it\lambda) m\ \in C$. We have the following well-know fact \cite{Bialynicki-Birula,Koras-86}.

\begin{proposition}
$C^-$ is a locally closed complex submanifold of $M$, and the projection $p_C:C^-\to C$ is an holomorphic map.
\end{proposition}

Suppose now that $\lambda\in\ugot^{-\sigma}$,  so $\lambda=i\beta$ with $\beta\in\pgot$. 
The complex submanifold $M^\lambda=M^\beta$ is then stable under the involution $\tau$. Suppose that 
$\Zcal^\beta=M^\beta\cap \Zcal$ is non-empty and let $\Ccal$ be a connected component of $\Zcal^\beta$. There exists a unique connected component 
$C\subset M^\lambda$ containing $\Ccal$ : the complex submanifold $C$ is then stable under the anti-holomorphic involution $\tau$, and 
$\Ccal$ is a connected component of $C^\tau$. We notice that for any $t\in\R$, and any $m\in M$ we have 
$\tau(\exp(-it\lambda) m)=\tau(\exp(t\beta) m)=\exp(t\beta) \tau(m)$. It shows that the complex Bialynicki-Birula's submanifolds $C^-$ is stable under $\tau$, hence 
$(C^-)^\tau=\{z\in M^\tau,\, \lim_{t\to\infty} \exp(t\beta)m\in C^\tau\}$ is a locally closed submanifold of $M^\tau$.

\begin{corollary}
$\Ccal^-=\{z\in \Zcal,\, \lim_{t\to\infty} \exp(t\beta)m\in \Ccal\}$ is a locally closed submanifold of $\Zcal$, and the projection $p_\Ccal:\Ccal^-\to \Ccal$ is a smooth map.
\end{corollary}

\subsection{Stratification in the complex setting}

 Let  
$$
f_\ugot:=\frac{1}{2}(\Phi_\ugot,\Phi_\ugot):M\longrightarrow \R
$$
denote the norm-square of the moment map. Notice that $f_\ugot$ is a proper function on $M$.

\begin{definition}
The Kirwan vector field $\kappa_\ugot$ is defined by the relation
$$
\kappa_\ugot(m)=\Phi_\ugot(m)\cdot m,\quad \forall m\in M.
$$
\end{definition}

We consider the gradient $\grad(f_\ugot)$ of the function $f_\ugot$ relatively to the Riemannian metric $\Omega(-,\J-)$. We recall the following well-known facts \cite{W11}.

\begin{proposition}
\begin{enumerate}
\item The gradient of $f_\ugot$ is $\grad(f_\ugot)=\J(\kappa_\ugot)$.
\item The set of critical points of the function $f_\ugot$ is $\critical(f_\ugot)=\{\kappa_\ugot=0\}$.
\item We have the decomposition
$\Phi_\ugot(\critical(f_\ugot))=\bigcup_{\lambda\in \Bcal_\Phi} U\lambda$ 
where the set $\Bcal_\ugot\subset \tgot^*_{+}$ is discrete. $\Bcal_\ugot$ is called the set of {\rm types} of $M$. 
\item We have  $\critical(f_\ugot)=\bigcup_{\lambda\in \Bcal_\ugot}\critical_\lambda$ 
where $\critical_\lambda=\critical(f_\ugot)\cap\Phi_\ugot^{-1}(U\lambda)$ is equal to $U(M^{\lambda}\cap \Phi_\ugot^{-1}(\lambda))$.
\end{enumerate}
\end{proposition}

Let $\varphi^t_\ugot:M\to M$ be the flow of $-\grad(f_\ugot)$; since $f_\ugot$ is proper, $\varphi^t_\ugot$ exists for all times
$t\in [0,\infty[$, and according to a result of Duistermaat \cite{Lerman05} we know that any trajectory of $\varphi^t_\ugot$ has a limit when $t\to\infty$. For any 
$m\in M$, let us denote $m_\infty:=\lim_{t\to\infty} \varphi^t_\ugot(m)$.

The construction of the Kirwan-Ness stratification goes as follows. For each $\lambda\in\Bcal_\ugot$, let $M_\lambda$ denote the set of points of $M$ flowing to $\critical_\lambda$ : 
$M_\lambda:=\{m\in M; m_\infty\in \critical_\lambda\}$. From its very definition, the set $M_\lambda$ is contained in $\{m\in M, f_\ugot(m)\geq \frac{1}{2}\|\lambda\|^2\}$. The Kirwan-Ness stratification is the decomposition \cite{Kirwan.84}, \cite{Ness84}: 
$$
M=\bigcup_{\lambda\in \Bcal_\ugot} M_\lambda.
$$

When $0$ belongs to the image of $\Phi_\ugot$, the strata $M_0$ corresponds to the dense open subset of analytical semi-stable points :
$M_0=\{m\in M; \overline{U_\C\, m}\cap \Phi_\ugot^{-1}(0)\neq\emptyset\}$.

Let us now explain the geometry of $M_\lambda$ for a non-zero type $\lambda$. Let $C_\lambda$ be the union of the connected components of $M^{\lambda}$ intersecting 
$\Phi_\ugot^{-1}(\lambda)$. Then $C_\lambda$ is a K\"{a}hler Hamiltonian $U_\lambda$-manifold\footnote{$U_\lambda$ is the stabilizer subgroup of $\lambda$.} with proper moment map $\Phi_\lambda:=\Phi_\ugot\vert_{C_\lambda}-\lambda$.

The Bialynicki-Birula's complex submanifold 
$$
C^-_\lambda:=\{m\in M,  \lim_{t\to\infty} \exp(-it\lambda)\cdot m\ \in C_\lambda\}
$$
corresponds to the set of points of $M$ flowing to $C_\lambda$ under the flow of $-\grad\langle\Phi_\ugot,\lambda\rangle$, as $t\to\infty$. The limit of the flow defines 
a projection $C^-_\lambda\to C_\lambda$. Notice that $C^-_\lambda$ is invariant under the action of the parabolic subgroup $P^\ugot_\lambda\subset U_\C$ : 
\begin{equation}\label{eq:P-gamma}
P_\lambda^\ugot=\{g\in U_\C, \lim_{t\to\infty}\exp(-it\lambda)g\exp(it\lambda)\ {\rm exists}\}.
\end{equation}
Consider now the Kirwan-Ness stratification of the K\"{a}hler Hamiltonian $U_\lambda$-manifold $C_\lambda$. Let $C_{\lambda,0}$ be the open strata of 
$C_\lambda$ corresponding to the $0$-type:
$$
C_{\lambda,0}=\{x\in C_{\lambda};\ \overline{(U_\lambda)_\C\, x}\cap \Phi^{-1}_\lambda(0)\neq \emptyset\}.
$$
Let $C^-_{\lambda,0}$ denote the inverse image of $C_{\lambda,0}$ in $C^-_\lambda$.

\begin{theorem}[Kirwan \cite{Kirwan.84}]\label{Kirwan-stratification}
Let $M$ be a K\"{a}hler Hamiltonian $U$-manifold with proper
moment map $\Phi_\ugot:M\to\ugot^*$. For each non zero type $\lambda$, $M_\lambda$ is a $U_\C$-invariant complex submanifold, and 
$U_\C\times_{P_\lambda^\ugot}C^-_{\lambda,0}\rightarrow M_\lambda$, $[g,z] \mapsto  g\cdot z$ 
is an isomorphism of complex  $U_\C$-manifolds.
\end{theorem}

Kirwan gave a proof when $M$ is a {\em compact} K\"{a}hler Hamiltonian $K$-manifold. When $M$ is non-compact but the moment map is proper, 
a proof is given in \cite{HSS08} (see also \cite{W11}). 

The following standard facts are going to be very important in the real setting. Let $\lambda_{\mini}$ be the orthogonal projection of 
 $0$ on the closed convex polytope $\Delta_\ugot(M)$. 

\begin{proposition}\label{prop:strate-ouverte-complexe}
\begin{enumerate}
\item[a)] $\lambda_{\mini}$ is the unique element of $\Bcal_\ugot$ with minimal norm.
\item[b)] $M_{\lambda_{\mini}}$ is a connected, open, dense and $U_\C$-invariant subset of $M$.
\item[c)] If $\lambda\neq \lambda_{\mini}$, then the strata $M_\lambda$ has an empty interior.
\end{enumerate}
\end{proposition}

\subsection{Stratification in the real setting : first step}

In this section we suppose that our K\"{a}hler Hamiltonian $U$-manifold $(M,\Omega,\J)$ admit an anti-holomorphic involution $\tau$ compatible with the 
complex-conjugate involution $\sigma$ on $U_\C$. The main purpose of this section is to show that the Kirwan-Ness stratification $M=\bigcup_{\lambda\in \Bcal_\ugot} M_\lambda$ 
induces a stratification on the real part $Z$ of $M$. In this section we suppose that the manifold $Z$ is {\em non-empty} : notice that the real dimension of any connected component of 
$Z$ is equal to the complex dimension of $M$.

Let us denote by $\sigma_+ :\tgot^*_+\to\tgot^*_+$ the involution of the Weyl chamber that is defined by the relations : $-\sigma(U\xi)=-U\sigma(\xi)=U\sigma_+(\xi)$ for any $\xi\in\tgot^*_+$.

Recall that $\varphi^t_\ugot:M\to M$ denotes the flow of $-\grad(f_\ugot)$. 
\begin{proposition}\label{prop:fundamental}
\begin{enumerate}
\item For any $(m,t)\in M\times \R_{\geq 0}$, we have $\tau(\varphi^t_\ugot(m))= \varphi^t_\ugot(\tau(m))$.
\item If $m_\infty=\lim_{t\to\infty}\varphi^t_\ugot(m)$, then $\tau(m_\infty)=(\tau(m))_\infty$ for all $m\in M$.
\item For any $\lambda\in\Bcal_\ugot$, we have 
\begin{itemize}
\item $\tau(M_\lambda)=M_{\sigma_+(\lambda)}$,
\item $M_\lambda\cap Z\neq\emptyset$ only if $\sigma(\lambda)=-\lambda$. 
\end{itemize}
\item If $\lambda\neq \lambda_{\mini}$, then the locally closed submanifold $M_\lambda\cap Z$ has an empty interior in $Z$.
\item $\sigma(\lambda_{\mini})=-\lambda_{\mini}$.
\item $M_{\lambda_{\mini}}\cap Z$ is a dense open subset of $Z$.
\item If $\Zcal$ is a connected component of $Z$, then $M_{\lambda_{\mini}}\cap \Zcal$ is a dense open subset of $\Zcal$, and the critical set $\critical_{\lambda_{\mini}}$ intersects $\Zcal$.
\end{enumerate}
\end{proposition}

{\em Proof :} The first point is a direct consequence of the fact that, for any $m\in M$, the tangent map $\T_m\tau :\T_m M\to\T_{\tau(m)}M$ sends $\J(\kappa_\ugot(m))$ to 
 $\J(\kappa_\ugot(\tau(m)))$. The second point follows from the first one.
 
 By definition $m\in M_\lambda \Leftrightarrow \Phi_\ugot(m_\infty)\in U\lambda$.  Then if $m\in M_\lambda$, we have
 $$
 \Phi_\ugot((\tau(m))_\infty)=\Phi_\ugot(\tau(m_\infty))=-\sigma(\Phi_\ugot(m_\infty))\in -U\sigma(\lambda)=U\sigma_+(\lambda).
 $$
Hence the identity $\tau(M_\lambda)=M_{\sigma_+(\lambda)}$ is proven. Let $m\in M_\lambda\cap Z$. Thanks to {\em 2.} we know that 
$m_\infty\in \Phi_\ugot^{-1}(U\lambda)\cap Z$ : it implies that $\Phi_\ugot(m_\infty)\in U\lambda\cap (\ugot^{-\sigma})^*$ and then $\sigma(\lambda)=-\lambda$ 
(see \S \ref{sec:symmetric-orbit}). The point {\em 3.} is settled.

Let $\lambda\neq \lambda_{\mini}$, and consider the decomposition of the (locally closed) complex submanifold $M_\lambda$ into connected components : 
$M_\lambda= \Ycal_1\cup\cdots\cup \Ycal_p$. If $M_\lambda\cap Z\neq \emptyset$, then $M_\lambda$ is stable under the involution $\tau$ : there exists a permutation 
$\epsilon \in\mathfrak{S}_p$ of order two such that $\tau(\Ycal_k)=\Ycal_{\epsilon(k)}$, and 
$$
M_\lambda\cap Z=\bigcup_{\epsilon(k)=k} \Ycal_k^{\tau}.
$$
Since $\dim_\C(\Ycal_k)<\dim_\C(M), \forall k$ (see point {\em c)} in Proposition \ref{prop:strate-ouverte-complexe}), we see that 
$\dim_\R(\Ycal_k^\tau)<\dim_\R(Z)$ for all $k$ such that $\epsilon(k)=k$. We can then conclude that the locally closed submanifold 
$M_\lambda\cap Z$ has an empty interior in $Z$. If $\Zcal$ is a connected component of $Z$, we prove similarly that $M_\lambda\cap \Zcal$ has an empty interior in $\Zcal$.

Consider the decomposition $Z= (M_{\lambda_{\mini}}\cap Z) \bigcup \bigcup_{\lambda\neq \lambda_{\mini}}(M_\lambda\cap Z)$, and 
the proper map $f_\ugot : M\to \R_{\geq 0}$ that we restrict to $Z$. If we denote $Z_{< R}=Z\cap \{f_\ugot < R\}$, we have  
$$
Z_{< R}= M_{\lambda_{\mini}}\cap Z_{< R} \bigcup \bigcup_{\lambda\neq \lambda_{\mini}, \|\lambda\|^2\leq 2R}M_\lambda\cap Z_{< R}.
$$
Since $M_\lambda\cap Z_{< R}$ has an empty interior in $Z_{<R}$, we see that $M_{\lambda_{\mini}}\cap Z_{< R}$ is dense in $Z_{< R}$ for any $R\geq 0$. Then
$M_{\lambda_{\mini}}\cap Z$ is dense in $Z$, and $M_{\lambda_{\mini}}\cap Z$ is open in $Z$ because $M_{\lambda_{\mini}}$ is open in $M$. In particular 
$M_{\lambda_{\mini}}\cap Z\neq \emptyset$ : hence $\sigma(\lambda_{\mini})=-\lambda_{\mini}$. We have proven the points {\em 5.} and {\em 6.}.

If $\Zcal$ is a connected component of $Z$, the decomposition $\Zcal= (M_{\lambda_{\mini}}\cap \Zcal) \bigcup 
\bigcup_{\lambda\neq \lambda_{\mini}}(M_\lambda\cap \Zcal)$ shows similarly that 
$M_{\lambda_{\mini}}\cap \Zcal$ is a dense open subset of $\Zcal$. If $x\in M_{\lambda_{\mini}}\cap \Zcal$ then $x_{\infty}\in \critical_{\lambda_{\mini}}\cap \Zcal$. 
We have proven point {\em 7.}. $\Box$

\subsection{Stratification in the real setting : second step}\label{sec:real-stratification}

We have a stratification 
$\Zcal= (M_{\lambda_{\mini}}\cap \Zcal)  \bigcup \bigcup_{\underset{\sigma(\lambda)=-\lambda}{\lambda\neq \lambda_{\mini}}}(M_\lambda\cap \Zcal)$
into locally closed submanifolds, that we are going to interpret in terms of the $G$-action on $\Zcal$ and the gradient map 
$\Phi_\pgot\vert_\Zcal:\Zcal\to\pgot^*$. Let $f_\pgot:=\tfrac{1}{2}(\Phi_\pgot,\Phi_\pgot):M\to\R$.

\begin{definition}
\begin{itemize}
\item If $\beta\in j^*(\Bcal_\ugot\cap (\ugot^{-\sigma})^*)$, we define $\Zcal_\beta= M_{(j^*)^{-1}(\beta)}\cap \Zcal$.
\item Let $\Bcal_\pgot=\{\beta\in j^*(\Bcal_\ugot\cap (\ugot^{-\sigma})^*), \Zcal_\beta\neq \emptyset\}\subset \agot^*_+$.
\item Let us denote $j^*(\lambda_{\mini})$ by $\beta_{\mini}$.
\end{itemize}
\end{definition}

\medskip

At this stage we have a stratification $\Zcal=\bigcup_{\beta\in\Bcal_\pgot}\Zcal_\beta$ into submanifolds. First let us check that $\Bcal_\pgot$ parametrizes the set of critical points of the function $f_\pgot\vert_\Zcal$. Recall that  point {\em 7.} in Proposition \ref{prop:fundamental} tell us that $0\in \Bcal_\pgot\Leftrightarrow (\Phi_\pgot\vert_\Zcal)^{-1}(0)\neq \emptyset\Leftrightarrow \Phi_\ugot^{-1}(0)\neq \emptyset$

\begin{lemma}
The set of critical points of the function $f_\pgot\vert_\Zcal$ admits the decomposition
$$
\critical(f_\pgot\vert_{\Zcal})=\bigcup_{\beta\in \Bcal_\pgot}K(\Zcal^{\beta}\cap \Phi_\pgot^{-1}(\beta)).
$$ 
\end{lemma}

{\em Proof :} We have $f_\ugot= \frac{1}{2}(\Phi_{\ugot^\sigma},\Phi_{\ugot^\sigma}) + f_\pgot$, and $\Phi_{\ugot^\sigma}:M\to(\ugot^\sigma)^*$ vanishes on $\Zcal$. Hence 
$\critical(f_\pgot\vert_\Zcal)=\critical(f_\ugot)\cap \Zcal$. Since $\critical(f_\ugot)=\bigcup_{\lambda\in \Bcal_\ugot}U(M^{\lambda}\cap \Phi_\ugot^{-1}(\lambda))$, we obtain
$\critical(f_\pgot\vert_\Zcal)=\bigcup_{\lambda\in \Bcal_\ugot} \Dcal_\lambda$ with $\Dcal_\lambda= U(M^{\lambda}\cap \Phi_\ugot^{-1}(\lambda))\bigcap \Zcal$. We see that 
$\Dcal_\lambda=\emptyset$ if $\lambda\notin \Bcal_\ugot\cap (\ugot^{-\sigma})^*$, otherwise $\Dcal_\lambda=K(\Zcal^{j^*(\lambda)}\cap \Phi_\pgot^{-1}(j^*(\lambda)))$ if 
$\lambda\in \Bcal_\ugot\cap (\ugot^{-\sigma})^*$. $\Box$

\medskip

Consider the case where $0\in\Bcal_\pgot$. 

\begin{lemma}\label{lem:Z-0}
If $0\in\Bcal_\pgot$ then $\Zcal_0=\left\{z\in \Zcal, \ \overline{G\, z}\cap \Phi_\pgot^{-1}(0)\neq\emptyset\right\}$. 
\end{lemma}

{\em Proof :} By definition,  $\Zcal_0=M_0\cap \Zcal= \left\{z\in \Zcal,\ z_\infty\in \Phi_\ugot^{-1}(0))\right\}$. 
Here $m_\infty=\underset{t\to\infty}{\lim}\varphi^t_\ugot(m)$, where $\varphi^t_\ugot:M\to M$ denotes the flow of $-\grad(f_\ugot)$. When $z\in \Zcal$, we know that 
$\varphi^t_\ugot(z)\in \Zcal, \forall t\geq 0$, and moreover the tangent vector 
$$
\frac{d}{dt}\varphi^t_\ugot(z)= -\grad(f_\ugot)(\varphi^t_\ugot(z))=-\grad(f_\pgot)(\varphi^t_\ugot(z))
$$
belongs to the subspace $\ggot\cdot\varphi^t_\ugot(z)=\{X\cdot\varphi^t_\ugot(z), X\in\ggot\}$. We see then that $\varphi^t_\ugot(z)\in Gz$ for all $t\geq 0$.
Since $\Zcal\cap \Phi_\ugot^{-1}(0)=\Zcal\cap \Phi_\pgot^{-1}(0)$, we have proved that $\Zcal_0\subset\left\{z\in \Zcal, \ \overline{G\, z}\cap \Phi_\pgot^{-1}(0)\neq\emptyset\right\}$. On the other hand, for any $z\in \Zcal$, 
$$
\overline{G\, z}\cap \Phi_\pgot^{-1}(0)\neq\emptyset\Longrightarrow \overline{U_\C\, z}\cap \Phi_\ugot^{-1}(0)\neq\emptyset \Longrightarrow z\in M_0.
$$
The proof is complete. $\Box$

Let $\beta\in\Bcal_\pgot$ be a non-zero element and let $\lambda=j^*(\beta)\in \Bcal_\ugot$. The strata $M_\lambda$ is isomorphic to $U_\C\times_{P_\lambda^\ugot}C^-_{\lambda,0}$.
We see that $\Ccal_\beta:= C_{\lambda}\cap \Zcal$ is the union of the connected components of $\Zcal^\beta$ intersecting $\Phi_\pgot^{-1}(\beta)$. Similarly,  
$\Ccal_\beta^-:= C_{\lambda}^-\cap \Zcal$ corresponds to {\em real} Bialynicki-Birula's submanifold 
$\{z\in \Zcal,  \lim_{t\to\infty} e^{t\beta} z\ \in \Ccal_\beta\}$. In the same way, $\Ccal_{\beta,0}:= C_{\lambda,0}\cap \Zcal$ is the open and dense subset of $\Ccal_{\beta}$ 
formed by the elements $z\in \Ccal_{\beta}$ such that $\overline{G_\beta\, z}\cap \Phi^{-1}_\pgot(\beta)\neq \emptyset$ (see Lemma \ref{lem:Z-0}). 
Finally, $\Ccal_{\beta,0}^-:=C_{\lambda,0}^-\cap \Zcal$ is the equal to the pullback of the open subset $\Ccal_{\beta,0}$ through the projection $\Ccal_{\beta}^-\to \Ccal_{\beta}$. 

The complex parabolic subgroup $P^\ugot_\lambda\subset U_\C$ is stable under the complex-conjugate involution $\sigma$ and the fixed point subgroup 
$(P^\ugot_\lambda)^\sigma$ is equal the real parabolic subgroup $P_\beta$ of $G$ defined by (\ref{eq:P-gamma-reel}). We can now conclude that the strata 
$\Zcal_\beta:=M_\lambda\cap \Zcal$ admits the following description.

\begin{proposition}\label{prop:description-strata}
Let $\beta\in\Bcal_\pgot$ be a non-zero element. The map $[g,z]\mapsto gz$ induces a diffeomorphism 
\begin{equation}\label{eq:decomposition-reelle}
G\times_{P_\beta}\Ccal_{\beta,0}^{-}\overset{\sim}{\longrightarrow} \Zcal_\beta.
\end{equation}
\end{proposition}

\medskip

The existence of a stratification $\Zcal=\bigcup_{\beta\in\Bcal_\pgot}\Zcal_\beta$ where each strata is described through the isomorphism 
(\ref{eq:decomposition-reelle}) was already obtained by P. Heinzner, G.W. Schwarz and H. St\"{o}tzel in a much more general setting \cite{HSS08}. 
Nevertheless we obtain in this particular framework a crucial information for the continuation: only one stratum has a non-empty interior, the open stratum $\Zcal_{\beta_{\mini}}$ attached to the element $\beta_{\mini}\in\Bcal_\pgot$ of minimal norm.

\subsection{Proof of the refined O'Shea-Sjamaar's theorem (in the K\"{a}hler case)}\label{sec:refinement}

Suppose that $Z\neq \emptyset$. Point $(6)$ of Proposition \ref{prop:fundamental} shows that $0\in \Delta_\ugot(M)$ if and only if $0\in \Delta_\pgot(Z)$. 
Let us explain how  O'Shea-Sjamaar's identify $\Delta_\ugot(M)\bigcap (\ugot^{-\sigma})^*\simeq\Delta_\pgot(Z)$ follows directly from this observation thanks to the shifting trick. 

The fact that $j^*(\Delta_\ugot(M)\bigcap (\ugot^{-\sigma})^*)$ contains $\Delta_\pgot(Z)$ is immediate. 

Take now $\lambda\in \Delta_\ugot(M)\bigcap (\ugot^{-\sigma})^*$, and consider the 
K\"{a}hler Hamiltonian $U$-manifold $N=M\times U(-\lambda)$. As $\sigma(\lambda)=-\lambda$, the map $(-\sigma)$ leaves $U(-\lambda)$ invariant and defines 
an anti-holomorphic involution $\tau_\lambda:U(-\lambda)\to U(-\lambda)$ by the relations $\tau_\lambda(-g\lambda)=-\sigma(g)\lambda$. 
By definition $Z_N=Z\times K(-\lambda)$ is the real part of $N$ and Point $(6)$ of Proposition \ref{prop:fundamental} shows that $0\in  \Delta_\pgot(Z_N)$ because $0\in \Delta_\ugot(N)$. 
We have proven that $j^{*}(\lambda)\in \Phi_\pgot(Z)\cap \agot^*_+$. The proof of the identity $\Delta_\ugot(M)\bigcap (\ugot^{-\sigma})^*\simeq\Delta_\pgot(Z)$ is complete.

The refinement of O'Shea-Sjamaar's theorem is a consequence Point $(7)$ of Proposition \ref{prop:fundamental} which says that if $0\in \Delta_\ugot(M)$ then the open subset 
$M_0$  of analytical semi-stable points intersects any connected component $\Zcal$ of $Z$. But if $x\in\Zcal\cap M_0$, then $x_\infty\in\Zcal\cap\Phi_\pgot^{-1}(0)$. 
In other words, $0\in \Delta_\ugot(M)$ if and only if $0\in \Delta_\pgot(\Zcal)$. Like before, the shifting trick shows that 
$j^*(\Delta_\ugot(M)\bigcap (\ugot^{-\sigma})^*)=\Delta_\pgot(\Zcal)$.

\section{Geometric properties}\label{sec:geometry}

\subsection{Symmetric coadjoint orbits}\label{sec:symmetric-orbit}

Let $\tilde{\lambda}\in\tgot^*$. The coadjoint  orbit $U\tilde{\lambda}$ admits a complex structure compatible with its symplectic structure which can be visualized through the isomorphisms $U\tilde{\lambda}\simeq U/U_{\tilde{\lambda}}\simeq U_\C/P^\ugot_{\tilde{\lambda}}$. Here  $P^\ugot_{\tilde{\lambda}}\subset U_\C$ is 
the complex parabolic subgroup defined by (\ref{eq:P-gamma}). Notice that the Lie algebra of $P^\ugot_{\tilde{\lambda}}$ is 
$$
{\rm Lie}(P^\ugot_{\tilde{\lambda}})=\tgot_\C\oplus \sum_{\underset{(\alpha,\tilde{\lambda})\leq 0}{\alpha\in\Rgot(\ugot,\tgot)}}(\ugot_\C)_{\alpha}.
$$
The conjugate linear involution $\sigma: \ugot_\C\to\ugot_\C$ leaves $\tgot_\C$ invariant and sends the weight space $(\ugot_\C)_{\alpha}$ to $(\ugot_\C)_{-\sigma(\alpha)}$. That permits to see that the image of the parabolic subgroup $P^\ugot_{\tilde{\lambda}}$ through $\sigma$ is the parabolic subgroup $P^\ugot_{-\sigma(\tilde{\lambda})}$.

Suppose now that coadjoint  orbit $U\tilde{\lambda}$ intersects $(\ugot^{-\sigma})^*$. If $\tilde{\lambda}$ is taken in the Weyl chamber, it is possible only if $\sigma(\tilde{\lambda})=-\tilde{\lambda}$, and in this case we have  $U\tilde{\lambda}\cap (\ugot^{-\sigma})^*=K\tilde{\lambda}$.  The map $\xi\in\ugot^*\mapsto -\sigma(\xi)\in\ugot^*$ leaves $U\tilde{\lambda}$ invariant and it defines an 
anti-holomorphic involution $\tau :U\tilde{\lambda}\to U\tilde{\lambda}$. The parabolic subgroup $P^\ugot_{\tilde{\lambda}}$ is stable under $\sigma$, so that the involution 
$\tau:U_\C/P^\ugot_{\tilde{\lambda}}\to U_\C/P^\ugot_{\tilde{\lambda}}$ can also be defined by $\tau([g])=([\sigma(g)])$.

The submanifold fixed by the involution $\tau :U{\tilde{\lambda}}\to U{\tilde{\lambda}}$ is $K{\tilde{\lambda}}$, and the isomorphism $U{\tilde{\lambda}}\simeq U_\C/P^\ugot_{\tilde{\lambda}}$ descends to an isomorphism $K\tilde{\lambda}\simeq G/ G\cap P^\ugot_{\tilde{\lambda}}$.

Let $j^*:\pgot^*\to (u^{-\sigma})^*$ be the $K$-equivariant isomorphism. Let $\lambda\in\agot^*$ such that $j^*(\lambda)=\tilde{\lambda}$. Then
$K\tilde{\lambda}\simeq K\lambda\simeq G/P_\lambda$ where $P_\lambda=G\cap P^\ugot_{\tilde{\lambda}}$ is the real parabolic subgroup of $G$ 
defined by (\ref{eq:P-gamma-reel}).

Let $\xi$ be an element belonging to the interior of the Weyl chamber $\agot^*_+$ : we have $(\alpha,\xi)>0$ for any restricted root $\alpha\in\Sigma^+$. In this case 
the orbit $K(-\xi)$ admits a natural identification with the quotient $G/Q$ since the (minimal) parabolic subgroup $Q$ coincides with $P_{-\xi}$ (see Definition \ref{def:parabolic-Q}).

\subsection{Construction of real Ressayre's pairs}\label{sec:construction-RP}

Let $\xi$ be an element in the interior of the Weyl chamber $\agot^*_{+}$ that does not belongs to the moment polytope $\Delta_\pgot(\Zcal)$. 
Let $\xi'$ be the orthogonal projection of $\xi$ on $\Delta_\pgot(\Zcal)$ and let $\gamma_\xi=\xi'-\xi\in\agot^*\simeq \agot$.

\begin{lemma}\label{lem:fibre}
The set $(\Phi_\pgot\vert_\Zcal)^{-1}(\xi')$ is contained in the submanidold $\Zcal^{\gamma_\xi}$.
\end{lemma}
{\em Proof :} Let us check that $\|\xi'-\xi\|^2$ is the minimal value of the function $\|\Phi_\pgot-\xi\|^2:\Zcal\to\R$. If $z\in \Zcal$ then 
$\Phi_\pgot(z)=k\eta$, with $k\in K$ and $\eta\in\Delta_\pgot(\Zcal)$ : $\|\Phi_\pgot(z)-\xi\|^2=\|\eta\|^2+\|\xi\|^2-2(k\eta,\xi)$. Since 
$(k\eta,\xi)\leq (\eta,\xi),\forall k\in K$, we get 
$$
\|\Phi_\pgot(z)-\xi\|^2\geq\|\eta\|^2+\|\xi\|^2-2(\eta,\xi)=\|\eta-\xi\|^2\geq \|\xi'-\xi\|^2.
$$
Then if $z\in \Phi_\pgot^{-1}(\xi')$, the differential of $\|\Phi_\pgot-\xi\|^2$ vanishes at $z$. But 
$$
d\|\Phi_\pgot-\xi\|^2\vert_z=2d\langle\Phi_\pgot,\gamma_\xi\rangle\vert_z=2(\gamma_\xi\cdot z,-).
$$ 
Thus $\gamma_\xi\cdot z=0$. $\Box$

\medskip

Let $\Ccal_{\gamma_\xi}$ be the union of the connected components of $\Zcal^{\gamma_\xi}$ intersecting $(\Phi_\pgot\vert_\Zcal)^{-1}(\xi')$. The next result is the main tool to exhibit 
real Ressayre's pairs on the $G$-manifold $\Zcal$.

\medskip

\begin{theorem}\label{theo:construction-RP}
The data $(\gamma_\xi, \Ccal_{\gamma_\xi})$ is a real Ressayre's pair on $\Zcal$.
\end{theorem}

The rest of this section is dedicated to the proof of Theorem \ref{theo:construction-RP}.

\medskip

Let $\tilde{\xi}$ be the element of $\tgot^*_+\cap(\ugot^{-\sigma})^*$ such that $\xi=j^*(\tilde{\xi})$. Since $\xi\notin\Delta_\pgot(\Zcal)$, we know that $\tilde{\xi}\notin \Delta_\ugot(M)$. 
Let $\tilde{\xi}'$ be the orthogonal projection of $\tilde{\xi}$ on $\Delta_\ugot(M)$. We work with the proper K\"{a}hler Hamiltonian $U$-manifold 
$$
N:=M\times U(-\tilde{\xi})\simeq  M\times U_\C/P^{\ugot}_{-\tilde{\xi}}.
$$
The submanifold $\Zcal_N:=\Zcal\times K(-\tilde{\xi})\simeq \Zcal\times G/Q$ is a connected component of its real part. 
Let $\Phi_\ugot^N:N\to \ugot^*$ be the moment map relative to the action of $U$ on $N$.

We following result precises Lemma \ref{lem:fibre}.

\begin{lemma}
\begin{itemize}
\item $\|g\tilde{\xi}'-\tilde{\xi}\|\geq \|\tilde{\xi}'-\tilde{\xi}\|, \, \forall g\in U$, and the equality holds if and only if $g\tilde{\xi}'\in U_{\tilde{\xi}}\cdot\tilde{\xi}'$.
\item The function $\|\Phi^N_\ugot\|:N\to\R$ reaches its minimum on $U(\Phi_\ugot^{-1}(\tilde{\xi}')\times\{-\tilde{\xi}\})$. 
\item We have $\sigma(\tilde{\xi}')=-\sigma(\tilde{\xi}')$. In other words, we have $\xi'=j^*(\tilde{\xi}')$ where $\xi'$ is the orthogonal projection of $\xi$ on $\Delta_\pgot(\Zcal)$.
\end{itemize}
\end{lemma}

{\em Proof :} The first point is a classical result of Hamiltonian geometry: we briefly recall the arguments. We have 
$\|g\tilde{\xi}'-\tilde{\xi}\|=\|\tilde{\xi}'\|^2+\|\tilde{\xi}\|^2-2\phi(g\tilde{\xi}')$ where $\phi$  is the
$\tilde{\xi}$-th component of the moment map on $U\tilde{\xi}'$. The function $\phi$ has a unique local maximum on the
coadjoint orbit $U\tilde{\xi}'$ which is reached on an orbit of the stabilizer subgroup $U_{\tilde{\xi}}$ (see \cite{Atiyah82,Guillemin-Sternberg82.bis}). Finally, it is not hard to check that 
the point $\tilde{\xi}'$ belongs to this orbit.

If $n=(m,k\tilde{\xi})\in N$, we write  $\Phi_\ugot(m)=g\eta$ with $\eta\in\Delta_\ugot(M)$. Then
$$
\|\Phi^N_\ugot(n)\|=\|g\eta-k\tilde{\xi}\|\geq \|\eta-\tilde{\xi}\|\geq \|\tilde{\xi}'-\tilde{\xi}\|
$$
and the equality $\|\Phi^N_\ugot(n)\|= \|\tilde{\xi}'-\tilde{\xi}\|$ holds if and only if $\eta=\tilde{\xi}'$ and $k^{-1}g\in U_{\tilde{\xi}}$. It follows that the critical set 
$\critical_{\lambda_{\mini}}=\{n,\, \|\Phi^N_\ugot(n)\|= \|\tilde{\xi}'-\tilde{\xi}\|\}$ is $U(\Phi_\ugot^{-1}(\tilde{\xi}')\times\{-\tilde{\xi}\})$. $\Box$

Thanks to point {\em 8)} of Proposition \ref{prop:fundamental}, we know that $\critical_{\lambda_{\mini}}\cap \Zcal_N$ is non-empty. 
If $g(m,-\tilde{\xi})\in \critical_{\lambda_{\mini}}\cap \Zcal_N$ then $gm\in \Zcal$ and $\Phi_\ugot(m)=\tilde{\xi}'$ which implies that 
$\Phi_\ugot(gm)$ belongs to $U\tilde{\xi}'\cap (\ugot^{-\sigma})^*$. It follows that $\sigma(\tilde{\xi}')=-\sigma(\tilde{\xi}')$ (see \S \ref{sec:symmetric-orbit}). 
$\Box$

\bigskip

Thanks to the previous lemma we know that the minimal type of $N=M\times U_\C/P^{\ugot}_{-\tilde{\xi}}$ is the non-zero element 
$\tilde{\gamma_\xi}:=\tilde{\xi}'-\tilde{\xi}\in (\tgot^{-\sigma})^*$: we 
denote by $\gamma_\xi:=j^*(\tilde{\gamma_\xi})=\xi'-\xi$ the corresponding element of $\agot^*\simeq \agot$. We will now use Proposition \ref{prop:description-strata} in order to describe the open and dense strata $(\Zcal_N)_{\gamma_\xi}$ of $\Zcal_N$.

\medskip

The critical set associated to the minimal type $\tilde{\gamma_\xi}$ is $(\Phi^N_\ugot)^{-1}(U\tilde{\gamma_\xi})=U(\Phi_\ugot^{-1}(\tilde{\xi}')\times\{-\tilde{\xi}\})$. Let 
$C_{N,\tilde{\gamma_\xi}}$ be the connected component of $N^{\tilde{\gamma_\xi}}$ containing $(\Phi^N_\ugot)^{-1}(U\tilde{\gamma_\xi})$ : we have
 $$
 C_{N,\tilde{\gamma_\xi}}=C_{\tilde{\gamma_\xi}}\times (U_{\tilde{\gamma_\xi}})_\C/((U_{\tilde{\gamma_\xi}})_\C\cap P^{\ugot}_{-\tilde{\xi}}),
 $$ 
 where $C_{\tilde{\gamma_\xi}}$ is the connected component of $M^{\tilde{\gamma_\xi}}$ containing $\Phi_\ugot^{-1}(\tilde{\xi}')$. The intersection 
 $C_{N,\tilde{\gamma_\xi}}\cap\Zcal_N$ decomposes as follow
$$
C_{N,\tilde{\gamma_\xi}}\cap\Zcal_N = \Ccal_{\gamma_\xi}\times G_{\gamma_\xi}/G_{\gamma_\xi}\cap Q, 
$$
where $\Ccal_{\gamma_\xi}$ is the union of the connected components of $\Zcal^{\gamma_\xi}$ intersecting $\Phi_\pgot^{-1}(\xi')$, and $G_{\gamma_\xi}\subset G$ 
is the subgroup that stabilizes $\gamma_\xi$. The real Bialynicki-Birula's submanifold  is then
$(C_{N,\tilde{\gamma_\xi}}\cap\Zcal_N)^-= \Ccal_{\gamma_\xi}^-\times P_{\gamma_\xi}/P_{\gamma_\xi}\cap Q$.

Let us consider the map 
\begin{equation}\label{map-p-gamma}
{\rm p}_{\gamma_\xi}: G\times_{P_{\gamma_\xi}}\left(\Ccal_{\gamma_\xi}^-\times P_{\gamma_\xi}/P_{\gamma_\xi}\cap Q\right)\longrightarrow \Zcal\times G/Q.
\end{equation}
defined by ${\rm p}_\gamma([g;z,[p]])=(gz,[gp])$. Proposition \ref{prop:description-strata} tells us that the image of ${\rm p}_{\gamma_\xi}$ is a dense 
$G$-invariant subset of $\Zcal\times G/Q$ and that there exists a dense $P_{\gamma_\xi}$-invariant subset $\Ucal$ of 
$\Ccal_{\gamma_\xi}^-\times P_{\gamma_\xi}/P_{\gamma_\xi}\cap Q$, 
intersecting $ \Ccal_{\gamma_\xi}\times G_{\gamma_\xi}/G_{\gamma_\xi}\cap Q$, and such that ${\rm p}_{\gamma_\xi}$ induces a diffeomorphism 
$G\times_{P_{\gamma_\xi}}\Ucal\simeq G\Ucal$.

\medskip

 In order finish the proof of Theorem \ref{theo:construction-RP}, we have to compare the map  ${\rm p}_{\gamma_\xi}$ defined by (\ref{map-p-gamma}) and the map 
${\rm q}_{\gamma_\xi}: Q\times_{Q\cap P_{\gamma_\xi}} \Ccal_{\gamma_\xi}^-\to \Zcal$.

Consider the  canonical isomorphism $\Ccal_{\gamma_\xi}^-\times P_{\gamma_\xi}/P_{\gamma_\xi}\cap Q \simeq P_{\gamma_\xi}\times_{P_{\gamma_\xi}\cap Q} \Ccal_{\gamma_\xi}^-$. 
More generaly, for any $P_{\gamma_\xi}$-invariant open subset 
$\Ucal\subset \Ccal_{\gamma_\xi}^-\times P_{\gamma_\xi}/P_{\gamma_\xi}\cap Q$ we have an isomorphism $\Ucal\simeq P_{\gamma_\xi}\times_{P_{\gamma_\xi}\cap Q} U$ where 
$U$ is the open $P_\gamma\cap Q$-invariant subset of $\Ccal^-_{\gamma_\xi}$ defined by the relation $U:=\{x\in \Ccal_{\gamma_\xi}^-; (x,[e])\in \Ucal\}$. 
Note that $\Ucal$ intersects $ \Ccal_{\gamma_\xi}\times G_{\gamma_\xi}/G_{\gamma_\xi}\cap Q$ if and only if $U$ intersects $\Ccal_{\gamma_\xi}$. 

Now, we notice that ${\rm q}_{\gamma_\xi}$ defines a diffeomorphism $Q\times_{Q\cap P_{\gamma_\xi}} U \simeq QU$ if and only if 
${\rm p}_{\gamma_\xi}$ defines a diffeomorphism $G\times_{P_{\gamma_\xi}} \Ucal \simeq G \Ucal$. That can be seen easily through the commutative diagram
$$
\xymatrixcolsep{5pc}\xymatrix{
G\times_{P_{\gamma_\xi}} \Ucal\ar[d]^{{\rm p}_{\gamma_\xi}} \ar[r]^{\sim} & G\times_{Q} (Q\times_{Q\cap P_{\gamma_\xi}} U)\ar[d]^{1\times{\rm q}_{\gamma_\xi}} \\
\Zcal\times G/Q \ar[r]^{\sim}          & G\times_Q \Zcal.
}
$$
The previous diagram shows also that ${\rm p}_{\gamma_\xi}$ has a dense image if and only if ${\rm q}_{\gamma_\xi}$ has a dense image. 
The proof of Theorem \ref{theo:construction-RP} 
is then complete. $\Box$

\subsection{Principal-cross-section Theorem}\label{sec:principal-cross-section}

We propose here a principal-cross-section Theorem in the spirit of \cite{L-M-T-W} that holds on the $K$-manifold $\Zcal$.  Since $\Delta_\pgot(\Zcal)$ is a closed convex polyhedral subset of $\agot^*_+$, 
there exists a unique open face $\sgot$ of $\agot^*_+$ such that  
\begin{itemize}
\item $\Delta_\pgot(\Zcal)\cap\sgot\neq \emptyset$,
\item $\Delta_\pgot(\Zcal)\subset\overline{\sgot}$.
\end{itemize}
Notice that $\Delta_\pgot(\Zcal)\cap\sgot$ is dense in $\Delta_\pgot(\Zcal)$. Recall that the linear isomorphism $j^*:(\tgot^*)^{-\sigma}\to \agot^*$ induces a bijection 
$\tgot^*_+\cap(\tgot^{-\sigma})^*\simeq \agot^*_+$.

\begin{lemma}
\begin{enumerate}
\item There exists a unique open face $\tilde{\sgot}$  of the Weyl chamber $\tgot^*_+$ such that $\tilde{\sgot}\cap(\tgot^{-\sigma})^*{\simeq}_{j^*} \sgot$.
\item The stabilizer subgroup $K_\xi$ does not depend of $\xi\in\sgot$ : it is denoted by $K_\sgot$.
\item The stabilizer subalgebra $\ggot_\xi$ does not depend $\xi\in\sgot$ : it is denoted by $\ggot_\sgot$. We have the decomposition : 
$\ggot_\sgot=\kgot_\sgot\oplus\pgot_\sgot$ where $\kgot_\sgot$ is the Lie algebra of $K_\sgot$.
\end{enumerate}
\end{lemma}

{\em Proof :} Any $\xi_o\in\tgot^*_+$ belongs to the open face $\tau(\xi_o)\subset \tgot^*_+$ defined as follows : $\xi\in \tau(\xi_o)$ if and only if  $\xi\in\tgot^*_+$ and 
$(\alpha,\xi)=0\Longleftrightarrow(\alpha,\xi_o)=0$ for any $\alpha\in \Rgot(\ugot,\tgot)$. Now, it is an easy matter to check that 
the face $\tau(\xi_o)$ does not depend on $\xi_o\in (j^*)^{-1}(\sgot)\subset\tgot^*_+\cap(\tgot^{-\sigma})^*$ : this face, denoted by $\tilde{\sgot}$, satisfies the relation 
$\tilde{\sgot}\cap(\tgot^{-\sigma})^* {\simeq}_{j^*} \sgot$.

All points in the open face $\tilde{\sgot}$ have the same connected stabilizer $U_{\tilde{\sgot}}$. If we take $\tilde{\xi}\in \tilde{\sgot}\cap(\tgot^{-\sigma})^*$, we see that 
$U_{\tilde{\xi}}=U_{\tilde{\sgot}}$ is invariant under $\sigma$ and that $K\cap U_{\tilde{\sgot}}=K\cap U_{\tilde{\xi}} $ is equal to $K_\xi$ where $\xi=j^*(\tilde{\xi})\in\sgot$.  
The second point is settled and the third point is leaved to the readers. $\Box$

\medskip

The following slice 
$$
\Tcal_\sgot=\{z\in\Zcal, \Phi_\pgot(z)\in \sgot\}.
$$
is the key object of our principal-cross-section Theorem.

\begin{theorem}\label{theo:principal-cross-section}
\begin{enumerate}
\item $\Tcal_\sgot$ is a $K_\sgot$-invariant submanifold of $\Zcal$.
\item The map $K\times_{K_\sgot}\Tcal_{\sgot}\to \Zcal$, $[k,y]\mapsto ky$ is a diffeomorphism onto a $K$-invariant open and dense subset of $\Zcal$.
\item For any $x\in \Tcal_\sgot$, we have $[\kgot_\sgot,\pgot_\sgot]\subset \pgot_x\subset \pgot_\sgot$.
\end{enumerate}
\end{theorem}

\medskip

{\em Proof :} We consider the following open subset of $\ugot_{\tilde{\sgot}}^*$ :
$$
\Vcal_{\tilde{\sgot}}:=U_{\tilde{\sgot}}\Big\{\xi\in\tgot^*_+, \ U_\xi\subset U_{\tilde{\sgot}}\Big\}.
$$

The pull-back $Y_{\tilde{\sgot}}=\Phi_\ugot^{-1}(\Vcal_{\tilde{\sgot}})$ is the symplectic cross-section at $\tilde{\sgot}$. It is an $U_{\tilde{\sgot}}$-invariant symplectic manifold of $M$ such that the map 
$U\times_{U_{\tilde{\sgot}}}Y_{\tilde{\sgot}}\longrightarrow M, [g,y]\mapsto gy$ defines a diffeomorphism onto the open and dense subset $U Y_{\tilde{\sgot}}$.

\begin{lemma}
\begin{enumerate}
\item $\Vcal_{\tilde{\sgot}}\subset \ugot_{\tilde{\sgot}}^*$ is invariant under the map $-\sigma$.
\item The submanifold $Y_{\tilde{\sgot}}$ is stable under the involution $\tau$.
\item The intersection $Y_{\tilde{\sgot}}\cap\Zcal$ is equal to $\Tcal_\sgot$.
\item The intersection $U Y_{\tilde{\sgot}}\cap\Zcal$ is equal to $K\Tcal_{\sgot}$.
\item $K\Tcal_{\sgot}$ is open and dense in $\Zcal$.
\end{enumerate}
\end{lemma}

{\em Proof of the lemma :} Recall that $K'$ is the connected component of the centraliser subgroup $Z_K(\agot)$. Let $w_0'$ be the longuest element in the Weyl group 
$W'=N_{K'}(T)/T$. Then the linear map $\sigma_+(\xi)=-w_0'\sigma(\xi)$ of $\tgot^*$ preserves the Weyl chamber $\tgot^*_+$, and for any $\xi\in \tgot^*_+$ ,we have 
$-\sigma(U\xi)=U(\sigma_+(\xi))$. Let $k'\in K'$ be a representant of $w_0'$. 

Notice that the subgroup $U_{\tilde{\sgot}}$ is stable under $\sigma$ since $-\sigma$ fixes the elements of $\tilde{\sgot}\cap(\tgot^{-\sigma})^*$. Take $\eta\in \Vcal_{\tilde{\sgot}}$: so 
$\eta=g\xi$ where $g\in U_{\tilde{\sgot}}$ and $\xi\in \tgot^*_+$ satisfies $U_\xi\subset U_{\tilde{\sgot}}$. We have 
$-\sigma(\eta)=\sigma(g)k'\sigma_+(\xi)$, where $\sigma(g)k'\in U_{\tilde{\sgot}}$ because $K'\subset U_{\tilde{\sgot}}$. Now we see that the stabilizer subgroup 
$U_{\sigma_+(\xi)}$ is equal to $Ad(k')\circ\sigma(U_\xi)$: as $U_\xi\subset U_{\tilde{\sgot}}$, we obtain $U_{\sigma_+(\xi)}\subset U_{\tilde{\sgot}}$. The first point is proved and the second point follows directly from the first.

The inclusion $\Tcal_\sgot\subset Y_{\tilde{\sgot}}\cap\Zcal$ is immediate. Let us check the reverse inclusion. Let $z\in Y_{\tilde{\sgot}}\cap\Zcal$ and $\eta=\Phi_\ugot(z)\in 
\Vcal_{\tilde{\sgot}}\cap (\ugot^{-\sigma})^*$. Taking the decomposition $\eta=g\xi$ as before, we see that $g\xi\in U\xi\cap(\ugot^{-\sigma})^*=K\xi$ : there exists $k\in K$ such that 
$g\xi=k\xi$ or in other words $g^{-1}k\in U_\xi\subset U_{\tilde{\sgot}}$. But $g\in U_{\tilde{\sgot}}$ and so $k\in U_{\tilde{\sgot}}\cap K=K_\sgot$. At this stage we know that 
$\Phi_\pgot(k^{-1}z)=j^*(\xi)\in \agot^*_+$. On one hand we know that $\xi\in\overline{\tilde{\sgot}}$ because $j^*(\xi)\in \Delta_\pgot(\Zcal)\subset \overline{\sgot}$. On the other hand we know $U_\xi\subset U_{\tilde{\sgot}}$. It shows that $\xi$ belongs to $\tilde{\sgot}$, and hence $k^{-1}z\in \Tcal_\sgot$. Since $k\in K_\sgot$ we can conclude that $z$ belongs to $\Tcal_\sgot$. 
The third point is settled. 

Thanks to the third point we know that the intersection $U Y_{\tilde{\sgot}}\cap\Zcal$ contains $K\Tcal_\sgot$. Let us prove that 
$U Y_{\tilde{\sgot}}\cap\Zcal\subset K\Tcal_\sgot$. Let $(z,y,u)\in\Zcal\times Y_{\tilde{\sgot}}\times U$ such that $z=uy$. We write $\Phi_\ugot(y)=g\xi$,  
where $g\in U_{\tilde{\sgot}}$ and $\xi\in \tgot^*_+$ satisfies $U_\xi\subset U_{\tilde{\sgot}}$. We see then that $\Phi_u(z)=ug\xi$ belongs to $U\xi\cap(\ugot^{-\sigma})^*=K\xi$ : 
there exists $k\in K$ such that $ug\xi=k\xi$, so $k^{-1}ug\in U_{\tilde{\sgot}}$. We have proved that there exists $g'\in U_{\tilde{\sgot}}$ such that $u=kg'$. The identity  
$z=kg'y$ shows then that $g'y=k^{-1}z\in Y_{\tilde{\sgot}}\cap\Zcal=\Tcal_\sgot$. We have proved that $z\in K\Tcal_\sgot$.

$K\Tcal_{\sgot}$ is open in $\Zcal$ since $U Y_{\tilde{\sgot}}$ is open in $M$.  Let $z\in \Zcal$, and assume that $\lambda=\Phi_\pgot(z)$ belongs to the Weyl chamber $\agot^*_+$. 
The {\em local convexity theorem} of Sjamaar (see \cite{sjamaar98}[Theorem 6.5] and \cite{OSS99}[Theorem 8.2]) tells us that for any $K$-invariant neighborhood $U\subset \Zcal$ of 
$z$, the local moment polytope $\Delta_\pgot(U)=\Phi_\pgot(U)\cap\agot^*_+$ is a neighborhood of $\lambda$ in $\Delta_\pgot(\Zcal)$. It implies that $U\cap \Tcal_{\sgot}\neq\emptyset$   
for any $K$-invariant neighborhood $U\subset \Zcal$ of $z$. The density of $K\Tcal_{\sgot}$ in $\Zcal$ is demonstrated. $\Box$

\medskip

We can now complete the proof of Theorem \ref{theo:principal-cross-section}.

The identity $Y_{\tilde{\sgot}}\cap\Zcal=\Tcal_\sgot$ shows that $\Tcal_\sgot$ corresponds to the union of the connected components of the submanifold $(Y_{\tilde{\sgot}})^\tau$ 
contained in $\Zcal$. Hence $\Tcal_\sgot$ is a submanifold of $\Zcal$.

The last point of the Lemma shows that the diffeomorphism $U\times_{U_{\tilde{\sgot}}}Y_{\tilde{\sgot}}\overset{\sim}{\longrightarrow} U Y_{\tilde{\sgot}}$ induces the diffeomorphism
$K\times_{K_{\sgot}}\Tcal_{\sgot}\overset{\sim}{\longrightarrow}K \Tcal_{\sgot}$.

Let us check the last point of Theorem \ref{theo:principal-cross-section}. For any $\beta\in \pgot$, the vector field $z\in\Zcal\mapsto \beta\cdot z$ is the gradient vector field of the function 
$\langle\Phi_\pgot,\beta\rangle :\Zcal\to \R$. Let $x\in\Tcal_\sgot$ :  then $\beta\in \pgot_x$ if and only if the differential $d\langle\Phi_\pgot,\beta\rangle\vert_x :\T_x\Zcal\to \R$ vanishes. Thanks to the second point of Theorem \ref{theo:principal-cross-section}, we know that 
$\T_x\Zcal=\T_x\Tcal_\sgot+\kgot\cdot x$. For any $X\in \kgot$, we have 
$$
d\langle\Phi_\pgot,\beta\rangle\vert_x(X\cdot x)=\langle\Phi_\pgot(x),[\beta,X]\rangle
$$
with $\Phi_\pgot(x)\in\sgot$. Thus $d\langle\Phi_\pgot,\beta\rangle\vert_x $ vanishes on $\kgot\cdot x$ if and only if $\beta\in\pgot_\sgot$. At this stage we have that 
$\pgot_x\subset \pgot_\sgot$. The function $\Phi_\pgot$, when restricted to the submanifold $\Tcal_\sgot$, takes value in $\sgot$. 
If we take $\beta\in[\kgot_\sgot,\pgot_\sgot]\subset \pgot_\sgot$, the function $\langle\Phi_\pgot,\beta\rangle$ is constant equal to zero on $\Tcal_\sgot$, thus
$d\langle\Phi_\pgot,\beta\rangle\vert_x$ vanishes on $\T_x\Tcal_\sgot$. We have checked that $[\kgot_\sgot,\pgot_\sgot]\subset\pgot_x$. 
$\Box$

\section{Proof of Theorem \ref{th:real-ressayre-pairs} }

Let $\Delta_\pgot(\Zcal)\subset\agot^*_{+}$ be the moment polytope of a connected component $\Zcal$ of the real part $Z$ (supposed non-empty) of a proper 
K\"{a}hler Hamiltonian $U$-manifold $(M,\Omega,\J)$.

 We define the following convex subsets of the chamber $\agot^*_{+}$. 

\begin{itemize}
\item $\Delta_{\infrp}$ is the set of points $\xi\in\agot^*_+$ satisfying  the inequalities $\langle \xi,\gamma\rangle\geq \langle \Phi_\pgot(\Ccal),\gamma\rangle$, for any 
{\bf infinitesimal  real Ressayre's pair} $(\gamma,\Ccal)$ of $\Zcal$.
\item $\Delta^{^{\reg}}_{\infrp}$ is the set of points $\xi\in\agot^*_+$ satisfying  the inequalities $\langle \xi,\gamma\rangle\geq \langle \Phi_\pgot(\Ccal),\gamma\rangle$, for any 
{\bf regular infinitesimal real Ressayre's pair} $(\gamma,\Ccal)$ of $\Zcal$.
\item $\Delta_{\rp}$ is the set of points $\xi\in\agot^*_+$ satisfying the inequalities $\langle \xi,\gamma\rangle\geq \langle \Phi_\pgot(\Ccal),\gamma\rangle$, 
for any {\bf real Ressayre's pair} $(\gamma,\Ccal)$ of $\Zcal$.
\item $\Delta^{^{\reg}}_{\rp}$ is the set of points $\xi\in\agot^*_+$ satisfying the inequalities $\langle \xi,\gamma\rangle\geq \langle \Phi_\pgot(\Ccal),\gamma\rangle$, 
for any {\bf regular real Ressayre's pair} $(\gamma,\Ccal)$ of $\Zcal$.
\end{itemize}

By definition, we have the commutative diagram, where all the maps are inclusions:
$$
\xymatrix{
\Delta_{\infrp} \ar@{^{(}->}[d] \ar@{^{(}->}[r] & \Delta^{^{\reg}}_{\infrp} \ar@{^{(}->}[d] &  \\
\Delta_{\rp} \ar@{^{(}->}[r]          & \Delta^{^{\reg}}_{\rp} \cdot }
$$

In \S \ref{sec:step-1} and \S \ref{sec:step-2}, we prove the inclusions $ \Delta_{\rp}\subset\Delta_\pgot(\Zcal)\subset \Delta_{\infrp}$. It follows then that 
$\Delta_\pgot(\Zcal)=\Delta_{\infrp}=\Delta_{\rp}$.

In \S \ref{sec:step-3}, we prove the inclusion $\Delta^{^{\reg}}_{\rp}\subset \Delta_\pgot(\Zcal)$, and since 
$\Delta_\pgot(\Zcal)=\Delta_{\infrp}\subset \Delta^{^{\reg}}_{\infrp} \subset\Delta_{\rp}^{^{\reg}}$, we get finally that 
$\Delta_\pgot(\Zcal)=\Delta^{^{\reg}}_{\rp}=\Delta^{^{\reg}}_{\infrp}$. The proof of Theorem \ref{th:real-ressayre-pairs} is complete.

\subsection{$\Delta_{\rp}\subset\Delta_\pgot(\Zcal)$}\label{sec:step-1}

Let $\xi_o\in\agot^*_{+}$ that does not belong to $\Delta_\pgot(\Zcal)$. The aim of this section is to prove that $\xi_o\notin \Delta_{\rp}$. In other words, 
we will show the existence of a real Ressayre's pair $(\gamma,\Ccal)$ of $\Zcal$ such that $\langle \xi_o,\gamma\rangle < \langle \Phi_\pgot(\Ccal),\gamma\rangle$.

Let $r>0$ be the distance between $\xi_o$ and $\Delta_\pgot(\Zcal)$, and let $\xi$ be an element in the interior of the Weyl chamber $\agot^*_{+}$ such as 
$\|\xi-\xi_o\|<\frac{r}{2}$ : so the distance between $\xi$ and  $\Delta_\pgot(\Zcal)$ is stricly larger than $\frac{r}{2}$.

Let $\xi'$ be the orthogonal projection of $\xi$ on $\Delta_\pgot(\Zcal)$ and let $\gamma=\xi'-\xi\in\agot^*\simeq \agot$. Let $\Ccal$ be the union of the connected components 
of $\Zcal^{\gamma}$ intersecting $\Phi_\pgot^{-1}(\xi')$. Thanks to Theorem \ref{theo:construction-RP}, we know that $(\gamma, \Ccal)$ is a real Ressayre's pair on $\Zcal$.

Using the fact that $\gamma=\xi'-\xi$, we compute 
\begin{eqnarray*}
\langle \xi_o,\gamma\rangle- \langle \Phi_\pgot(\Ccal),\gamma\rangle&=&\langle \xi_o,\gamma\rangle- \langle \xi',\gamma\rangle\\
&=&\langle \xi_o-\xi,\gamma\rangle- \|\gamma\|^2\\
&\leq &- \|\gamma\|(\|\gamma\|- \|\xi_o-\xi\|)\\
&<&0\, .
\end{eqnarray*}
The last inequality comes from the fact that $\|\xi_o-\xi\|<\frac{r}{2}$ and that $\|\gamma\|>\frac{r}{2}$ since $\|\gamma\|$ represents the distance between $\xi$ and $\Delta_\pgot(\Zcal)$.

\subsection{$\Delta_\pgot(\Zcal)\subset \Delta_{\infrp}$}\label{sec:step-2}

Let $\xi\in\agot^*_{+}$ belonging to $\Delta_\pgot(\Zcal)$. The aim of this section is to prove that $\xi\in \Delta_{\infrp}$. 
Let $\tilde{\xi}$ be the corresponding element in $\tgot^*_+\cap(\ugot^{-\sigma})^*$: in other words $\xi=j^*(\tilde{\xi})$. 
Consider the K\"{a}hler Hamiltonian  $U$-manifold $N:=M\times U(-\tilde{\xi})$ and the connected component $\Zcal_N:=\Zcal\times K(-\tilde{\xi})$ of its real part. 

Let $\Phi_\pgot^N: \Zcal_N\to\pgot^*$ be the gradient map. We know that the strata 
$$
(\Zcal_N)_0=\left\{n\in \Zcal_N,\, \overline{G\, n}\cap (\Phi_\pgot^N)^{-1}(0)\neq \emptyset\right\}
$$  
is a $G$-invariant dense open subset of $\Zcal_N$.

Let $(\gamma, \Ccal)$ be an infinitesimal real Ressayre's pair on $\Zcal$, and let $\Ccal_N:= \Ccal\times (G_\gamma/P_{-\xi}\cap G_\gamma)$ be the 
corresponding connected component of $\Zcal_N^\gamma$. Notice that $\Ccal_N$ is invariant under the action of the stabilizer subgroup $G_\gamma$.

Let $\Ccal^-_N:=\{n\in \Zcal_N,  \lim_{t\to\infty} e^{t\gamma} n\ \in \Ccal_N\}$ be the real Bialynicki-Birula's submanifold.

\begin{lemma}\label{lem:KCinterior}
\begin{enumerate}
\item The set $G\Ccal^{-}_N$ has a non-empty interior in $\Zcal_N$. 
\item $\Ccal^-_N\cap (\Zcal_N)_0\neq \emptyset$.
\end{enumerate}
\end{lemma}
{\em Proof :}
The second point follows from the first one since $(\Zcal_N)_0$ is a dense $G$-invariant subset of $\Zcal_N$.

Let $x\in \Ccal$ so that $\ngot^{\gamma>0}\simeq (\T_x \Zcal)^{\gamma>0}$. The point $n=(x,[e])\in  \Ccal_N$ belongs to the interior of $G\Ccal^-_N$ if we show that 
$\ggot\cdot n + \T_n  \Ccal_N^-=\T_n \Zcal_N$. Since $\T_n \Ccal_N^-=(\T_n \Zcal_N)^{\gamma\leq 0}$ it is sufficient to check that 
$(\T_n \Zcal_N)^{\gamma>0}\subset \ggot\cdot n$. We have the decomposition $(\T_n \Zcal_N)^{\gamma>0}=(\T_x \Zcal)^{\gamma>0}\oplus \ggot^{\gamma>0}\cdot [e]$. 
So for any $v\in (\T_n \Zcal_N)^{\gamma>0}$, there exists $X\in\ggot^{\gamma>0}$ so that $v-X\cdot(x,[e])\in (\T_x  \Zcal_N)^{\gamma>0}$. 
But  $\ngot^{\gamma>0}\simeq (\T_x  \Zcal)^{\gamma>0}$, so there exists $Y\in \ngot^{\gamma>0}$ such as $v-X\cdot(x,[e])=Y\cdot x$. 
The Lie algebra $\ngot$ is contained in the Lie algebra of the parabolic subgroup $P_{-\xi}$: 
hence $Y\cdot [e]=0$. Finally we have proved that $v=(X+Y)\cdot(x,[e])\in  \ggot\cdot n$. $\Box$

\medskip

The proof of the inclusion $\Delta_\pgot(\Zcal)\subset \Delta_{\infrp}$ follows from the next result.

\begin{lemma}Let $(\gamma, \Ccal)$ be an infinitesimal real Ressayre's pair on $\Zcal$.
Then, the following inequality holds  $\langle \xi,\gamma\rangle \geq \langle \Phi_\pgot(\Ccal),\gamma\rangle$ for any $\xi\in\agot^*_{+}$ belonging to $\Delta_\pgot(\Zcal)$.
\end{lemma}

{\em Proof :} Let $n\in \Ccal^-_N\cap (\Zcal_N)_0$, and let $n_\gamma\in\Ccal_N$ be the limit $\lim_{t\to\infty}e^{t\gamma} n$.

Let $P_\gamma\subset G$ be the parabolic subgroup associated to $\gamma$ (see (\ref{eq:P-gamma})). Since $G=KP_\gamma$, the fact 
that $n\in (\Zcal_N)_0$ means that $\overline{P_\gamma n}\cap(\Phi^N_\pgot)^{-1}(0)\neq \emptyset$. In other words,
$\min_{z\in P_\gamma n}\|\Phi^N_\pgot(z)\|=0$ which implies $0\geq \min_{z\in P_\gamma n}\langle\Phi^N_\pgot(z),\gamma\rangle$.

Consider now the function $t\geq 0\mapsto \langle \Phi^N_\pgot(e^{t\gamma} z),\gamma\rangle$ attached to $z\in P_\gamma n$.  Since 
$\frac{d}{dt}\langle\Phi^N_\pgot(e^{t\gamma}z),\gamma\rangle= -\|\gamma_{\Zcal_N}\|^2(e^{t\gamma}z)\leq 0$, 
we have 
\begin{equation}\label{eq:N-ss-gamma}
\langle\Phi^N_\pgot(z),\gamma\rangle \geq \langle\Phi^N_\pgot(e^{t\gamma}z),\gamma\rangle,\qquad \forall t\geq 0.
\end{equation} 

Let's take $p\in P_\gamma$ and $z=pn$. Then, the limit $\lim_{t\to\infty}e^{t\gamma} z$ is equal to $g_\gamma n_\gamma$ where 
$g_\gamma=\lim_{t\to\infty}e^{t\gamma} p \, e^{-t\gamma}\in G_\gamma$. If we take the limit  in (\ref{eq:N-ss-gamma}) as $t\to\infty$, we get
$$
\langle\Phi^N_\pgot(z),\gamma\rangle\geq \langle \Phi^N_\pgot(g_\gamma n_\gamma),\gamma\rangle= \langle \Phi^N_\pgot(\Ccal_N),\gamma\rangle=\langle \Phi_\pgot(\Ccal),\gamma\rangle-\langle \xi,\gamma\rangle,\quad \forall z\in P_\gamma n.
$$ 
We obtain finally $0\geq \min_{z\in P_\gamma n}\langle\Phi^N_\pgot(z),\gamma\rangle\geq \langle \Phi_\pgot(\Ccal),\gamma\rangle-\langle \xi,\gamma\rangle$.
$\Box$

\subsection{$\Delta^{^{\reg}}_{\rp}\subset \Delta_\pgot(\Zcal)$}\label{sec:step-3}

The aim of this section is the proof of the following 
\begin{theorem}
Let $\xi\in\agot^*_{+}$ satisfying the inequalities $\langle \xi,\gamma\rangle\geq \langle \Phi_\pgot(\Ccal),\gamma\rangle$, for any regular real Ressayre's pair $(\gamma,\Ccal)$ of $\Zcal$. 
Then $\xi\in\Delta_\pgot(\Zcal)$.
\end{theorem}

Our arguments go as follows: we will show that there exists a collection $(\gamma_i, \Ccal_i)_{i\in I}$ of regular real Ressayre's pairs
for which we have 
$$
\bigcap_{i\in I}\Big\{\xi\in\agot^*_{+}, \langle \xi,\gamma_i\rangle\geq \langle \Phi(\Ccal_i),\gamma_i\rangle\Big\}=\Delta_\pgot(\Zcal).
$$
The set $I$ will be finite when $\Zcal$ is compact.

Before starting the description of the collection $(\gamma_i, \Ccal_i)_{i\in I}$, recall the following fact concerning admissible elements (see \S 4.2 in \cite{pep-RP}). 
\begin{remark}\label{lem:admissible}
Let $\gamma\in\agot$ be a rational element such that  $K\,\Zcal^\gamma=\Zcal$. We have $\dim_\pgot(\Zcal^\gamma)=\dim_\pgot(\Zcal)$, thus $\gamma$ is admissible.
\end{remark}

\subsubsection{The real Ressayre's pair $(\gamma_\sgot, \Ccal_\sgot)$}\label{sec:RP-gamma-s}

In this section, we  show that a real Ressayre's pair describes the fact that $\Delta_\pgot(\Zcal)$ is contained in $\bar{\sgot}$. 
Let us denote by $h_\alpha\in \agot$ the coroot associated to a root $\alpha\in \Sigma$ :  $h_\alpha$ is the rational element of $[\ggot_\alpha,\ggot_{-\alpha}]\cap\agot$ satisfying 
$\langle \alpha, X\rangle=(h_\alpha,X)_b$, $\forall X\in\agot$ (see the Appendix).

 Let $\Sigma^+\subset \Sigma$ be the set of positive roots associated to the choice of the Weyl chamber $\agot^*_+$ : $\xi\in \agot^*_+$ if and only if 
 $\langle \xi, h_\alpha\rangle\geq 0$ for every $\alpha\in \Sigma^+$. Let $\Sigma^+_\sgot\subset \Sigma^+$ be the set of positive roots orthogonal to $\sgot$ :
 $\alpha\in \Sigma^+_\sgot$ if $\langle \xi, h_\alpha\rangle = 0$ for every $\xi\in\sgot$.
 
\begin{definition}
Consider the following rational vector 
$$
\gamma_\sgot:= -\sum_{\alpha\in\Sigma^+_\sgot}h_\alpha \ \in \agot.
$$
\end{definition}

\begin{lemma}The element $\gamma_\sgot$ satisfies the following properties:
\begin{itemize}
\item $\langle\xi,\gamma_\sgot\rangle\leq 0$ for any $\xi\in\agot^*_{\geq 0}$.
\item For any $\xi\in\agot^*_{\geq 0}$, $\langle\xi,\gamma_\sgot\rangle=0$ if and only if $\xi\in \bar{\sgot}$.
\item $\gamma_\sgot$  acts trivially on the principal-cross-section\footnote{See \S \ref{sec:principal-cross-section}.} $\Tcal_\sgot$.
\item $\langle\alpha,\gamma_\sgot\rangle< 0$ for any $\alpha\in\Sigma^+_\sgot$.
\end{itemize}
\end{lemma}

{\em Proof :} The first two points follow from the fact that $\langle\xi,h_\alpha\rangle\geq 0$ for any $\xi\in\agot^*_{\geq 0}$ and any positive roots $\alpha$. 
The third point is due to the fact that  $\gamma_\sgot\in [\kgot_\sgot,\pgot_\sgot]$. For the last point, see the Appendix. $\Box$

Let $\Ccal_\sgot$ be the union of the connected components of $\Zcal^{\gamma_\sgot}$ intersecting $\Tcal_\sgot$. We start with the following basic result.

\begin{lemma}\label{lem:gamma-s}
\begin{itemize}
\item $\gamma_\sgot$ is an admissible element.
\item For any $\xi\in\agot^*_+$, the inequality $\langle \xi,\gamma_\sgot\rangle\geq \langle \Phi_\pgot(\Ccal_\sgot),\gamma_\sgot\rangle$
is equivalent to $\xi\in\bar{\sgot}$.
\end{itemize}
\end{lemma}

{\em Proof :} 
We know that $\Tcal_\sgot\subset\Ccal_\sgot$ and that $K\Tcal_\sgot$ is dense in $\Zcal$. It follows that  $K \Ccal_\sgot=\Zcal$, hence 
$\gamma_\sgot$ is an admissible element (see Remark \ref{lem:admissible}). 
Now we consider the inequality $\langle \xi,\gamma_\sgot\rangle\geq \langle \Phi_\pgot(\Ccal_\sgot),\gamma_\sgot\rangle$ for an element $\xi\in\agot^*_{+}$. First we notice that
$\langle \Phi_\pgot(\Ccal_\sgot),\gamma_\sgot\rangle=0$, and the first two points of the previous Lemma tell us that $\langle \xi,\gamma_\sgot\rangle\geq 0$ is equivalent to 
$\xi\in\bar{\sgot}$.
$\Box$

\begin{proposition}\label{prop:gamma-tau}
$(\gamma_\sgot, \Ccal_\sgot)$ is a regular real Ressayre's pair on $\Zcal$.
\end{proposition}

{\em Proof :}  Using the identification $\agot\simeq \agot^*$, we view $\gamma_\sgot$ as a rational element of $\sgot^*$ orthogonal to $\sgot$. 
Take $\xi'\in\Delta_\pgot(\Zcal)\cap\sgot$ contained in the image $\Phi_\pgot(\Ccal_\sgot)$ and consider the elements $\xi(n):=\xi'-\frac{1}{n}\gamma_\sgot$ for $n\geq 1$. 
We notice that for $n$ large enough
\begin{enumerate}
\item $\xi(n)$ is a regular element of the Weyl chamber $\agot^*_+$,
\item $\xi(n)\notin \Delta_\pgot(\Zcal)$,
\item $\xi'$ is the orthogonal projection of $\xi(n)$ on $\Delta_\pgot(\Zcal)$.
\end{enumerate}
So we can exploit the results of \S \ref{sec:construction-RP} with the elements $\xi(n)$ for $n>>1$. Proposition  
\ref{theo:construction-RP} and Lemma \ref{lem:gamma-s} tell us that $(\gamma_\sgot, \Ccal_\sgot)$ is a regular real Ressayre's pair.
$\Box$

\subsubsection{The real Ressayre's pairs $(\gamma^\pm_l, \Ccal^\pm_l)$}

Let $\R\sgot\subset\agot^*$ be the rational vector subspace generated by the face $\sgot$. The closed convex polytope $\Delta_\pgot(\Zcal)$ generate an affine subspace 
$\Pi$ of $\R\sgot$. In this section we show that a finite family of real Ressayre's pairs describe the fact that $\Delta_\pgot(\Zcal)$, viewed as a subset of $\R\sgot$, is contained in 
the affine subspace $\Pi$.

In this section we will use the identifications $\pgot^*\simeq\pgot$, and $\agot\simeq\agot^*$ given by the invariant scalar product $(-,-)_b$ (see the Appendix).

We start with the orthogonal decompositions $\pgot=\agot\oplus\qgot$ and $\agot=\R\sgot\oplus\R\sgot^\perp$. It follows that $\pgot_\sgot=\agot\oplus\qgot_\sgot$. 
Let $\overrightarrow{\Pi}^\perp$ be the orthogonal of $\overrightarrow{\Pi}$ in $\R\sgot$.

For any $x\in \Zcal$, we define $\agot_x=\{X\in\agot, X\cdot x=0\}$.

\begin{lemma}\label{lem:agot-pi}
\begin{enumerate}
\item For any $x\in \Tcal_\sgot$, we have $\pgot_x=\agot_x\oplus\qgot_\sgot$ with $\overrightarrow{\Pi}^\perp\oplus\R\sgot^\perp\subset \agot_x$.
\item The equality $\overrightarrow{\Pi}^\perp\oplus\R\sgot^\perp= \agot_x$ holds on an open subset of $\Tcal_\sgot$, thus 
\begin{equation}\label{eq:stabilizer-generique}
\dim_\pgot(\Zcal)= \dim(\overrightarrow{\Pi}^\perp)+\dim(\R\sgot^\perp)+\dim(\qgot_\sgot).
\end{equation}
\item The subspaces $\overrightarrow{\Pi}^\perp$ and $\overrightarrow{\Pi}$ are rational.
\end{enumerate}
\end{lemma}

{\em Proof :} Let us come back to the arguments used in the proof of Theorem \ref{theo:principal-cross-section}. We know that $\pgot_x\subset\pgot_\sgot=\agot\oplus\qgot_\sgot$ 
for any $x\in  \Tcal_\sgot$, thus $\pgot_x=\agot_x\oplus\qgot_\sgot$. An element $\beta$ belongs to $\agot_x$ if and only if the differential 
$d\langle\Phi_\pgot,\beta\rangle\vert_x :\T_x\Tcal_\sgot\to \R$ vanishes. We see that for any $\beta\in \overrightarrow{\Pi}^\perp\oplus\R\sgot^\perp$ the map 
$\langle\Phi_\pgot,\beta\rangle :\Tcal_\sgot\to \R$ is constant : we obtain then that $\overrightarrow{\Pi}^\perp\oplus\R\sgot^\perp\subset \agot_x$, $\forall x\in \Tcal_\sgot$.

Since the image of the map  $\Phi_\pgot:\Tcal_\sgot\to \Pi$ is open in $\Pi$, there exists $x_o\in \Tcal_\sgot$ such that $d\Phi_\pgot\vert_{x_o}:\T_{x_o}\Tcal_\sgot\to \overrightarrow{\Pi}$ 
is surjective. It follows that $\overrightarrow{\Pi}^\perp\oplus\R\sgot^\perp= \agot_{x_o}$. The second point is proved.

The vector subspaces $\R\sgot\subset\agot$ and $\R\sgot^\perp\subset\agot$  are rational. Let $x_o\in\Tcal_\sgot$ such that  $\agot_{x_o}=\overrightarrow{\Pi}^\perp\oplus\R\sgot^\perp$. The Lie algebra 
$\tgot_{x_o}$ is a rational subspace of $\tgot$, thus $\agot_{x_o}=i(\tgot_{x_o}\cap\tgot^{-\sigma})$ is a rational subspace of $\agot$. If follows that 
$\overrightarrow{\Pi}^\perp=\agot_{x_o}\cap \R\sgot$ is a rational subspace of $\agot$. $\Box$

\medskip

Let $(\eta_l)_{l\in L}$ be a rational basis of $\overrightarrow{\Pi}^\perp$. We consider then the rational elements
$$
\gamma^\pm_l:= \pm \eta_l +\gamma_\sgot\ \in\ \overrightarrow{\Pi}^\perp\oplus\R\sgot^\perp,\quad l\in L.
$$
Thanks to Lemma \ref{lem:agot-pi}, we know that $\Tcal_\sgot\subset \Zcal^{\gamma^\pm_l}$, $\forall l\in L$. Let $\xi'_o\in \Delta_\pgot(\Zcal)\cap \sgot$. For any $l\in L$, we denote by 
$\Ccal^\pm_l$ the union of connected components of $\Zcal^{\gamma^\pm_l}$ intersecting $\Phi_\pgot^{-1}(\xi'_o)\subset\Tcal_\sgot$.

\begin{lemma}\label{lem:gamma-S}
\begin{itemize}
\item Any $\gamma^\pm_l$ is an admissible element.
\item The set of elements $\xi\in\R\sgot$ satisfying the inequalities
\begin{equation}\label{eq:gamma-S}
\langle \xi,\gamma_l^\pm\rangle\geq \langle \Phi_\pgot(\Ccal^\pm_l),\gamma^\pm_l\rangle, \quad \forall l\in L
\end{equation}
corresponds to the affine subspace $\Pi$.
\end{itemize}
\end{lemma}

{\em Proof :} Like in the proof of Lemma \ref{lem:gamma-s}, we see that $K\cdot \Zcal^{\gamma^\pm_l}=\Zcal$ since $\Tcal_\sgot\subset\Zcal^{\gamma^\pm_l}$ and $K\Tcal_\sgot$ is dense in $\Zcal$. It follows that $\gamma_\sgot$ is an admissible element.

The element $\xi'_o$ belongs to $\Pi$. Since $\langle\xi,\gamma_\sgot\rangle=0, \forall \xi\in\R\sgot$, the inequalities (\ref{eq:gamma-S}) are equivalent to
$\pm\langle \xi-\xi'_o,\eta_l\rangle\geq 0,  \forall l\in L$: in other words $\xi-\xi'_o\in \overrightarrow{\Pi}$, so
$\xi\in\Pi$.
$\Box$

\medskip

\begin{proposition}\label{prop:gamma-S}
For any $l\in L$, $(\gamma^\pm_l, \Ccal^\pm_l)$ is a regular real Ressayre's pair of $\Zcal$.
\end{proposition}

{\em Proof :} The proof follows the lines of the proof of Proposition \ref{prop:gamma-tau}.  The element $\xi'_o$ is 
contained in the image $\Phi_\pgot(\Ccal^\pm_l)$. We consider the elements $\xi_l^\pm(n):=\xi'_o-\frac{1}{n}\gamma_l^\pm$ for $n\geq 1$. 
We notice that for $n$ large enough
\begin{enumerate}
\item $\xi_l^\pm(n)$ is a regular element of the chamber $\agot^*_+$,
\item $\xi_l^\pm(n)\notin \Delta_\pgot(\Zcal)$,
\item $\xi'_o$ is the orthogonal projection of $\xi_l^\pm(n)$ on $\Delta_\pgot(\Zcal)$.
\end{enumerate}
So we can exploit the results of \S \ref{sec:construction-RP} with the elements $a_l^\pm(n)$ for $n>>1$. Proposition  
\ref{theo:construction-RP} and Lemma \ref{lem:gamma-S} tell us that $(\gamma^\pm_l, \Ccal^\pm_l)$  is a regular real Ressayre's pair.
$\Box$

\subsubsection{The real Ressayre's pairs $(\gamma_F, \Ccal_F)$}

In this section, we show that the polytope $\Delta_\pgot(\Zcal)$, viewed as a subset of the affine subspace $\Pi$, is the intersection of 
the cone $\Pi\cap\agot^*_{+}$, with a collection of half spaces of $\Pi$ parametrized by a family of real Ressayre's pairs.

\begin{definition}
An open facet $F$ of $\Delta_\pgot(\Zcal)$ is called {\em non trivial} if $F\subset\sgot$. We denote by $\Fcal(\Zcal)$ the set of non trivial open facets of $\Delta_\pgot(\Zcal)$.
\end{definition}

Let $F\in \Fcal(\Zcal)$. There exists\footnote{Here we see $\overrightarrow{\Pi}$ as a subspace of $\agot$ through the identification $\agot^*\simeq \agot$.} $\eta_F\in\overrightarrow{\Pi}$ such as the affine space generated by $F$ is 
$\Pi_F=\{\xi\in \Pi,\langle\xi,\eta_F\rangle=\langle\xi_F,\eta_F\rangle\}$ for any $\xi_F\in F$. The vector $\eta_F$ is chosen so that 
$\Delta_\pgot(\Zcal)$ is contained in the half space $\{\xi\in \Pi,\langle\xi,\eta_F\rangle\geq\langle\xi_F,\eta_F\rangle\}$.

By definition of the set $\Fcal(\Zcal)$, we have the following description of the polytope $\Delta_\pgot(\Zcal)$ :
\begin{equation}\label{eq:description-delta-K-M}
\Delta_\pgot(\Zcal)=  \bigcap_{F\in \Fcal(\Zcal)}\Big\{\xi\in \Pi,\langle\xi,\eta_F\rangle\geq\langle\xi_F,\eta_F\rangle\Big\}\hspace{2mm}\bigcap\hspace{2mm}\agot^*_{+}.
\end{equation}

The function $\langle\Phi_\pgot,\eta_F\rangle: \Tcal_\sgot\to\R$ is locally constant on $\Tcal_\sgot^{\eta_F}$ and it takes it's minimal value on 
$\Phi^{-1}_\pgot(F)\subset \Tcal_\sgot$, thus $\Phi^{-1}_\pgot(F)$ is an open subset of the submanifold $\Tcal_\sgot^{\eta_F}$. The map 
$\Phi_\pgot:\Phi^{-1}_\pgot(F)\to F$ is surjective, so it admits a regular value $\xi_F'\in F$.  For any $x\in \Phi^{-1}_\pgot(\xi'_F)$, 
the tangent map $\T_{x}\Phi_\pgot: \T_{x}\Tcal_\sgot^{\eta_F}\to \overrightarrow{F}$ is surjective, thus 
\begin{equation}\label{eq:stabilizer-F}
\agot_{x}= \R\eta_F\oplus \overrightarrow{\Pi}^\perp\oplus\R\sgot^\perp, \quad \forall x\in \Phi^{-1}_\pgot(\xi'_F).
\end{equation}
Since the vector subspaces $\agot_{x}$, $\overrightarrow{\Pi}^\perp$, $\R\sgot^\perp$ are rational, the vector $\eta_F$ can be taken rational.

We consider now the rational elements
$$
\gamma_F=\eta_F+\gamma_\sgot,\quad F\in\Fcal(\Zcal).
$$

Let $\Ccal_F$ be the union of the connected components of $\Zcal^{\gamma_F}$ intersecting $\Phi^{-1}_\pgot(\xi'_F)$.  The identities 
(\ref{eq:stabilizer-generique}) and (\ref{eq:stabilizer-F}) show that $\dim_\pgot(\Ccal_F)=\dim_\pgot(\Zcal)+1$.

\begin{lemma}\label{lem:gamma-F}
\begin{itemize}
\item $\gamma_F$ is an admissible element, for any $F\in \Fcal(\Zcal)$.
\item The set of elements $\xi\in\Pi\cap \tgot^*_{\geq 0}$ satisfying the inequalities
\begin{equation}\label{eq:gamma-F}
\langle \xi,\gamma_F\rangle\geq \langle \Phi_\pgot(\Ccal_F),\gamma_F\rangle, \quad \forall F\in \Fcal(\Phi)
\end{equation}
corresponds to $\Delta_\pgot(\Zcal)$.
\end{itemize}
\end{lemma}

{\em Proof :}
The first point is due to the fact that $\eta_F$ is rational and that $\dim_\pgot(\Ccal_F)=\dim_\pgot(\Zcal)+1$. 
The last assertion is a consequence of (\ref{eq:description-delta-K-M}). Notice that the relations 
$\langle \xi,\gamma_F\rangle\geq \langle \Phi_\pgot(\Ccal_F),\gamma_F\rangle$ and 
$\langle \xi,\eta_F\rangle\geq \langle \Phi_\pgot(\Ccal_F),\eta_F\rangle$ are equivalent for any $\xi\in\Pi$. $\Box$

\begin{proposition}\label{prop:gamma-F}
For any $F\in \Fcal(\Zcal)$, the couple $(\gamma_F, \Ccal_F)$ is a real Ressayre's pair on $\Zcal$.
\end{proposition}

{\em Proof :}
The proof follows the lines of the proof of Proposition \ref{prop:gamma-tau}.  
Consider the elements $\xi_F(n):=\xi'_F -\frac{1}{n}\gamma_F$ for $n\geq 1$. 
We notice that for $n$ is large enough
\begin{enumerate}
\item $\xi_F(n)$ is a regular element of the Weyl chamber,
\item $\xi_F(n)\notin  \Delta_\pgot(\Zcal)$,
\item $\xi'_F$ is the orthogonal projection of $\xi_F(n)$ on $ \Delta_\pgot(\Zcal)$.
\end{enumerate}
So we can exploit the results of \S \ref{sec:construction-RP} with the elements $\xi_F(n)$ for $n>>1$. Proposition  
\ref{theo:construction-RP} and Lemma \ref{lem:gamma-F} tell us that $(\gamma_F, \Ccal_F)$ is a regular real Ressayre's pair on $\Zcal$.
$\Box$

\section{Appendix}

In this paper, $U_\C$ denote a connected complex reductive Lie group with maximal compact subgroup $U$. Let $\sigma$ be a complex-conjugate involution on 
$U_\C$ leaving the sugbroup $U$ invariant. Let $T$ be a maximal torus of $U$ that is stable under $\sigma$ : we suppose that the vector subspace 
$\tgot^{-\sigma}=\{X\in\tgot,\sigma(X)=-X\}$ has a maximal dimension.

Take a faithful representation $\rho :U\to U(V)$ where $V$ is an Hermitian vector space. It extends to a morphism $\rho :U_\C\to {\rm GL}(V)$, so that $\rho(U_\C)$ is a closed subgroup stable under the Cartan involution $\Theta:{\rm GL}(V)\to{\rm GL}(V)$.

 We consider the real bilinear map 
\begin{equation}\label{eq:bilinear}
b:\ugot_\C\times \ugot_\C\longrightarrow \R
\end{equation}
defined by $b(X,Y)= {\rm Re}\Big({\rm Tr}(d\rho(X)d\rho(Y)) +{\rm Tr}(d\rho(\sigma(X))d\rho(\sigma(Y)))\Big)$. We see then that 
\begin{itemize}
\item $b(\sigma(X),\sigma(Y))=b(X,Y)$, $\forall X,Y\in\ugot_\C$, 
\item $b(uX,uY)=b(X,Y)$, $\forall u\in U_\C$, $\forall X,Y\in\ugot_\C$, 
\item $b(iX,iY)=-b(X,Y)$, if $X,Y\in\ugot_\C$,
\item $b(X,Y)\in \Z$, if $X,Y\in\wedge:=\frac{1}{2\pi}\ker(\exp:\tgot\to T)$.
\end{itemize}

In the article, we work with the following $U$-invariant scalar product on $\ugot_\C$ :
$$
(X,Y)_b= -b(X,\Theta(Y))
$$
If $V\subset \ugot_\C$ is a real vector subspace, we have an isomorphism $\xi\in V^*=\hom(V,\R)\to \xi^\flat \in V$ defined by the relation 
$$
\langle\xi, X\rangle=(\xi^\flat,X)_b,\qquad \forall (\xi,X)\in V^*\times V.
$$

Let $G$ be the real reductive group equal to the connected component of $(U_\C)^\sigma$: its maximal compact subgroup is $K=G\cap U$. At the level of Lie algebras we have 
$\ggot=\kgot\oplus \pgot$ where $\kgot=\ugot^{\sigma}$ and $\pgot=i\ugot^{-\sigma}$. The bilinear map (\ref{eq:bilinear}) defines a $G$-invariant bilinear map 
$b:\ggot\times\ggot\to\R$ that is positive definite on $\pgot$ and negative definite on $\kgot$.  The commutativity of the following diagram is frequently used in the 
body of the present paper :
$$
\xymatrixcolsep{5pc}\xymatrix{
(\ugot^{-\sigma})^*\ar[d]^{\flat} \ar[r]^{j^*} & \pgot^*\ar[d]^{\flat} \\
\ugot^{-\sigma}          & \pgot\ar[l]^{j} .
}
$$

We consider the restricted root system $\Sigma=\Rgot(\ggot,\agot)$ associated to a maximal abelian subspace $\agot\subset\pgot$ and a system of positive roots $\Sigma^+$ associated 
to a choice of Weyl chamber $\agot^*_+$.  Let $\sgot$ be a face of $\agot^*_+$: the set of positive roots orthogonal to $\sgot$ is denoted by $\Sigma^+_\sgot$. We consider 
$$
\rho_\sgot:=\sum_{\alpha\in\Sigma^+_\sgot}\alpha\quad \in\agot^*
$$
and the corresponding dual element $(\rho_\sgot)^\flat=-\gamma_\sgot\in \agot$ (see \S \ref{sec:RP-gamma-s} where the element $\gamma_\sgot$ is used).

\begin{lemma}\label{lem:appendix}
We have $\langle\alpha,\gamma_\sgot\rangle < 0$ for any $\alpha\in\Sigma^+_\sgot$.
\end{lemma}

{\em Proof :} Let $\alpha_o$ be a simple root of $\Sigma^+_\sgot$. Let ${\rm s}_{\alpha_o}:\agot^*\to\agot^*$ be the associated orthogonal symmetry. 
For any $\beta\in \Sigma^+_\sgot$, standard computations give that ${\rm s}_{\alpha_o}(\beta)\in \Sigma^+_\sgot$ if $\beta$ is not proportional to $\alpha$, and 
${\rm s}_{\alpha_o}(\beta)=-\beta$ if $\beta$ is proportional to $\alpha$. It shows that 
${\rm s}_{\alpha_o}(\rho_\sgot)=\rho_\sgot-2N\alpha_o$ for some integers $N\geq 1$, thus $(\alpha_o,\rho_\sgot)_b>0$. This implies that 
$(\alpha,\rho_\sgot)_b>0$ for any $\alpha\in \Sigma^+_\sgot$. $\Box$

\end{document}